\documentclass[letterpaper, preprint, paper,11pt]{AAS}	

\usepackage{bm}
\usepackage{amsmath}
\usepackage{amssymb}
\usepackage{subfigure}
\usepackage[colorlinks=true, pdfstartview=FitV, linkcolor=black, citecolor= black, urlcolor= black]{hyperref}
\usepackage{overcite}
\usepackage{footnpag}			      	
\usepackage{algorithm,algcompatible}
\algnewcommand\INPUT{\item[\textbf{Input:}]}%
\algnewcommand\OUTPUT{\item[\textbf{Output:}]}%

\PaperNumber{24-165}

\begin{document}

\title{Analysis and Design of Satellite Constellation Spare Strategy Using Markov Chain}

\author{Seungyeop Han\thanks{PhD Student, Daniel Guggenheim School of Aerospace Engineering, Georgia Institute of Technology, Atlanta, GA,30332},  
Takumi Noro\thanks{Researcher, Advanced Technology R\&D Center, Mitsubishi Electric Corporation, Amagasaki 661-8861, Japan.},
\ and Koki Ho\thanks{Dutton-Ducoffe Professor, Daniel Guggenheim School of Aerospace Engineering, Georgia Institute of Technology, Atlanta,GA, 30332}
}

\maketitle{}

\begin{abstract}
This paper introduces the analysis and design method of an optimal spare management policy using Markov chain for a large-scale satellite constellation. We propose an analysis methodology of spare strategy using a multi-echelon $(r,q)$ inventory control model with Markov chain, and review two different spare strategies: direct resupply, which inserts spares directly into the constellation orbit using launch vehicles; and indirect resupply, which places spares into parking orbits before transferring them to the constellation orbit. Furthermore, we propose an optimization formulation utilizing the results of the proposed analysis method, and an optimal solution is found using a genetic algorithm.
\end{abstract}

\section{Introduction}
Satellite constellations are groups of satellites working together over time as a system. Unlike single satellites, constellations can provide global coverage, with satellites in many different orbital planes. Due to their unprecedented usefulness, mega-constellations, which consist of more than a couple of thousand satellites, are being actively developed these days. Two representative projects of the New Space paradigm are OneWeb and Starlink. For example, OneWeb plans to launch 6,362 small satellites in low Earth orbit (LEO) to provide broadband services \cite{oneweb}, whereas nearly 12,000 interlinked broadband Internet satellites are planned to be deployed for Starlink, with a possible later extension to 42,000 \cite{starlink1, starlink2}.

However, the number of operating satellites will decrease over time due to failures, leading to a degradation in the constellation's performance. Therefore, an effective population management strategy is needed to maintain the constellation's designed performance. Two representative strategies for constellation management are the spare strategy and on-orbit servicing. The spare strategy involves inserting a new spare satellite to replace a failed one \cite{jakob2019optimal}, while on-orbit servicing entails sending a servicer to repair a failed satellite \cite{luu2022orbit}. Each method has its strengths and weaknesses; generally, the spare strategy is suitable for mega-constellations with small, low-cost satellites, while on-orbit servicing is more appropriate for smaller constellations comprising larger, high-cost satellites.

Among the two methods, this paper focuses on the spare strategy to address the optimal management policy for mega-constellation. Previous work \cite{cornara1999satellite} proposed different spare strategies and their approximate replacement times, and research\cite{jakob2019optimal} analyzed these strategies using a $(r,q)$ inventory control model. The study\cite{jakob2019optimal} examined the costs associated with maintaining operational satellite constellations, including supplying spares from in-plane redundant spares, spares in lower altitude parking orbits, and spares launched from the ground as needed. Additionally, the work\cite{novak2023analysis} adapted the result to analyze the resilience of constellation against solar weather radiation effects.

We propose a novel methodology for modeling and analyzing various spare strategies. Our approach utilizes Markov processes to model state transition behavior and analyze long-term behavior using stationary state distribution. This approach can handle the finite resource model, state-dependent failure model, and spare transfer constraints, which were ignored in previous work. In addition, the assumption regarding the spare demand probability distribution, which arises from satellite failures, is relaxed, resulting in a more accurate solution. Lastly, the proposed method computes the entire state probability distribution, allowing for various cost analyses.

With this method, we provide solutions for two main strategies: direct resupply, which inserts spares directly into the constellation orbit using launch vehicles; and indirect resupply, which inserts spares into parking orbits before transferring them to the constellation orbit. Using the results from our analysis method, we propose a simple optimization formulation that minimizes the maintenance cost of the system while satisfying the satellite population requirements.

The remainder of the paper is structured as follows. First, we provide the preliminaries of this work. Based on these preliminaries, we introduce the analysis method for the direct and indirect resupply strategies, respectively. Then, the analysis method is validated by comparing results from the Monte Carlo simulation, and the design optimization, which uses the analysis method, is demonstrated. Finally, we conclude the paper with a plan for future work.

\section{Preliminaries and Assumptions} \label{sec2}
\subsection{Orbital Mechanic}
\subsubsection{Constellation Model}
This research focuses on large-scale constellations in low Earth orbit (LEO), particularly the well-known Walker Delta pattern constellation\cite{walker1984satellite}. This constellation comprises multiple circular in-plane orbits with identical inclination angles. Additionally, the right ascension of the ascending node (RAAN) angle $\Omega$ of the in-plane orbits (referred to as the constellation orbits) is evenly distributed. In this configuration, there are a total of $N_\text{planes}$ in-plane orbits, with nominal $N_\text{sats}$ satellites allocated to each in-plane orbit.

\subsubsection{Parking Orbits}
The parking orbits are temporary orbits where the batch of spares is first inserted for the indirect resupply strategy. There are a total of $N_\text{park}$ parking orbits, and they are evenly distributed in RAAN. The parking orbits are assumed to have the same inclination angle as the in-plane orbits but have different semi-major axes.

\subsubsection{Orbital Maneuver}
This research assumes that all orbit transfers are conducted using coplanar Hohmann transfers, and out-of-plane maneuvers are excluded due to their inefficiency in terms of delta-v. Given the assumption of a LEO constellation, the time of flight for a Hohmann transfer is at most a couple of hours, which is negligible compared to the time horizon of the management policy. Therefore, it is assumed that the maneuver is instantaneous for simplicity.

\subsubsection{RAAN Drift}
The LEO experiences non-negligible perturbations due to the oblateness of the Earth, resulting in the secular drift of the orbital plane. In this problem, only the RAAN value will drift over time as a function of:
\begin{equation} \label{eq:rann_drift}
    \frac{d\Omega}{dt} = -\frac{3nR_\oplus^2 J_2}{2a^2} \cos i
\end{equation}
where $n$ is the mean motion, $a$ is the semi-major axis, $i$ is the inclination of the satellite orbit, $R_\oplus$ is the Earth's radius, and $J_2$ is the second zonal harmonic coefficient of the Earth\cite{prussing1993orbital}.

From this equation, it follows that an orbit with the same inclination but a different semi-major axis compared to the in-plane orbit will have a different RAAN drift value. Specifically, if the parking orbit and in-plane orbit have different semi-major axis, then the parking orbit will periodically align with all the constellation orbits, enabling efficient orbital transfer from the parking to the in-plane orbit.

\subsection{Spare Management Policy}
For this research, two different types of resupply methods are considered: direct and indirect resupply strategies. The first involves using a small LV to replenish the in-plane orbit directly. The second method utilizes a large LV to resupply the in-plane orbit indirectly via the parking orbit. The advantage of the second method over the first is that it can insert spares with significantly lower launch costs due to batch discounts. However, it will require considerably more time for replenishment due to the small RAAN drift rate. Given this trade-off relationship, it is crucial to compare the two methods across various case scenarios to determine the optimal spare management policy.

\begin{figure}[!h]
    \centering
    \includegraphics[width=.45\textwidth]{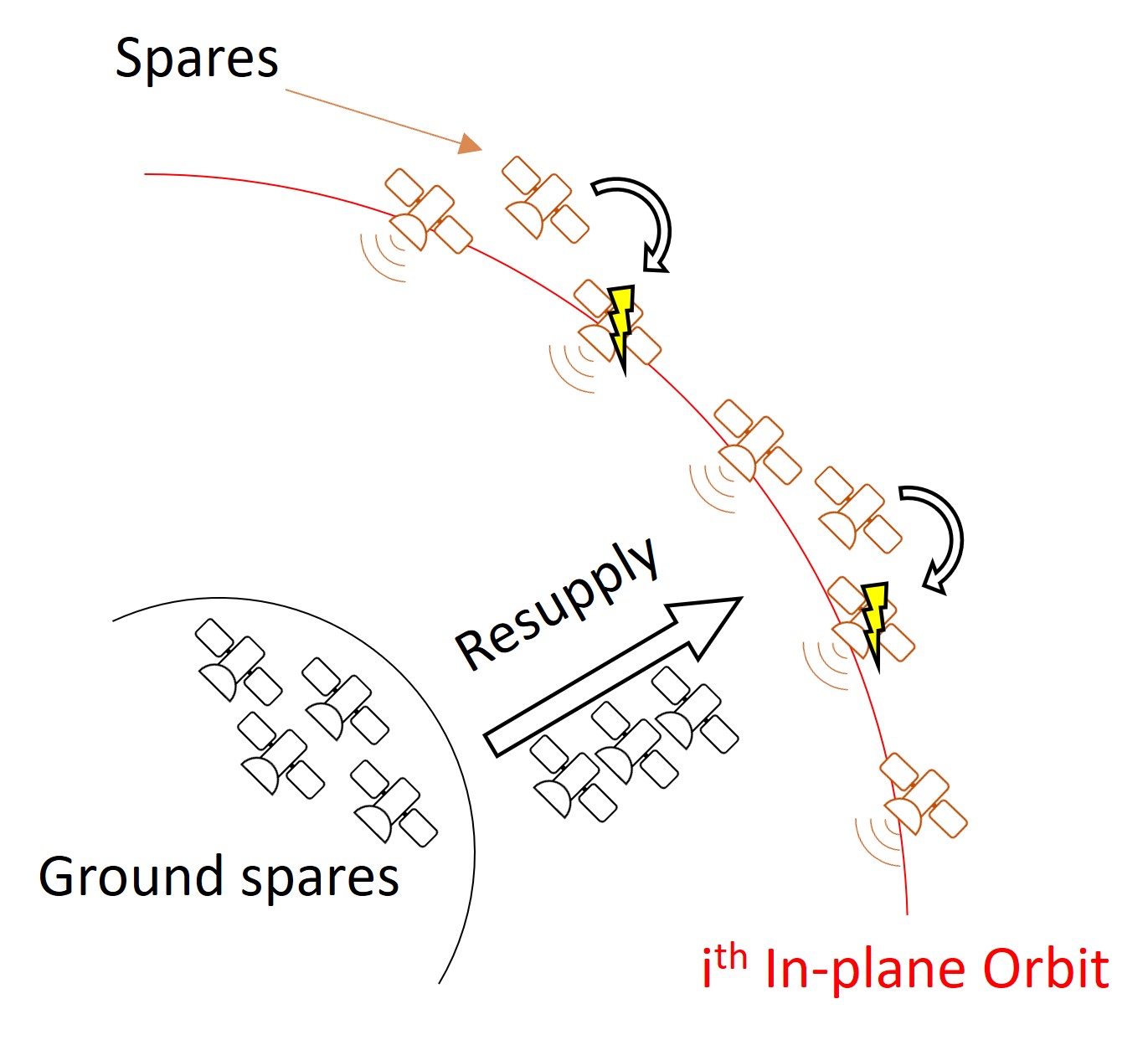}
    \includegraphics[width=.45\textwidth]{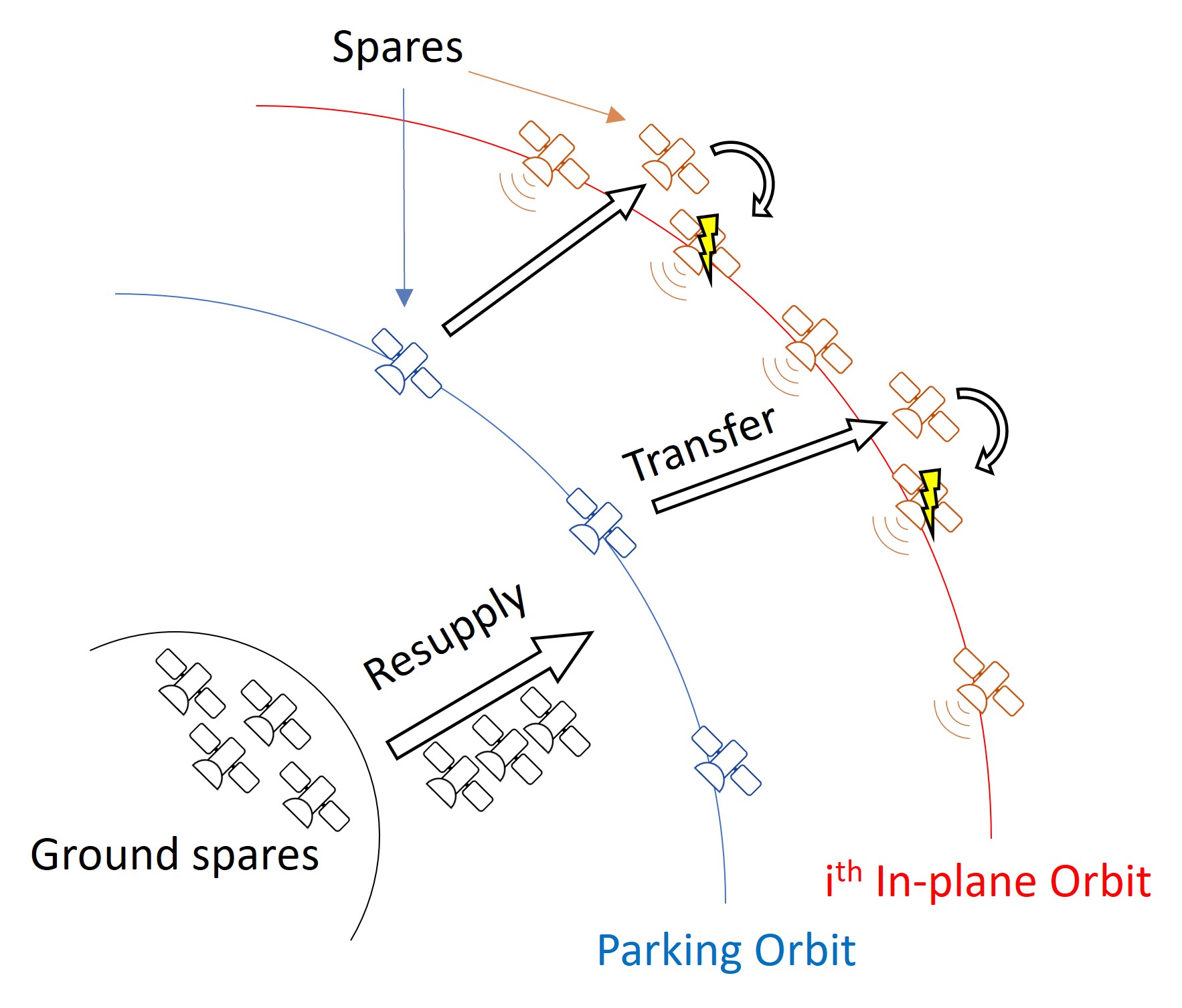}
    \caption{Two Different Replenishment Strategies (a) Direct Resupply (b) Indirect Resupply}
    \label{fig_strategy}
\end{figure}

The two different replenishment strategies are illustrated in Figure~\ref{fig_strategy}. In the direct resupply strategy, if a failure occurs, an existing in-plane orbit spare immediately replaces the failed one. Simultaneously, if the number of in-plane orbit spares falls below a certain threshold, a ground resupply order is placed. The resupply is then delivered after the lead time of a small-sized LV. For the indirect resupply strategy, the approach to replacing the failed satellite with an in-plane orbit spare remains the same. However, in this strategy, the ground resupply is first inserted into the parking orbit to replenish the parking orbit spares using a large-sized LV, possibly with a longer lead time. These parking orbit spares are then distributed from the parking orbit to the in-plane orbit as needed.

\subsection{Inventory Control Model}
To define the spare management policy rigorously, we introduce the following inventory control model.

\subsubsection{$(r,q)$ policy} \label{sec:single_rq}
\begin{figure}[h]
    \centering
    \includegraphics[width=0.5\textwidth]{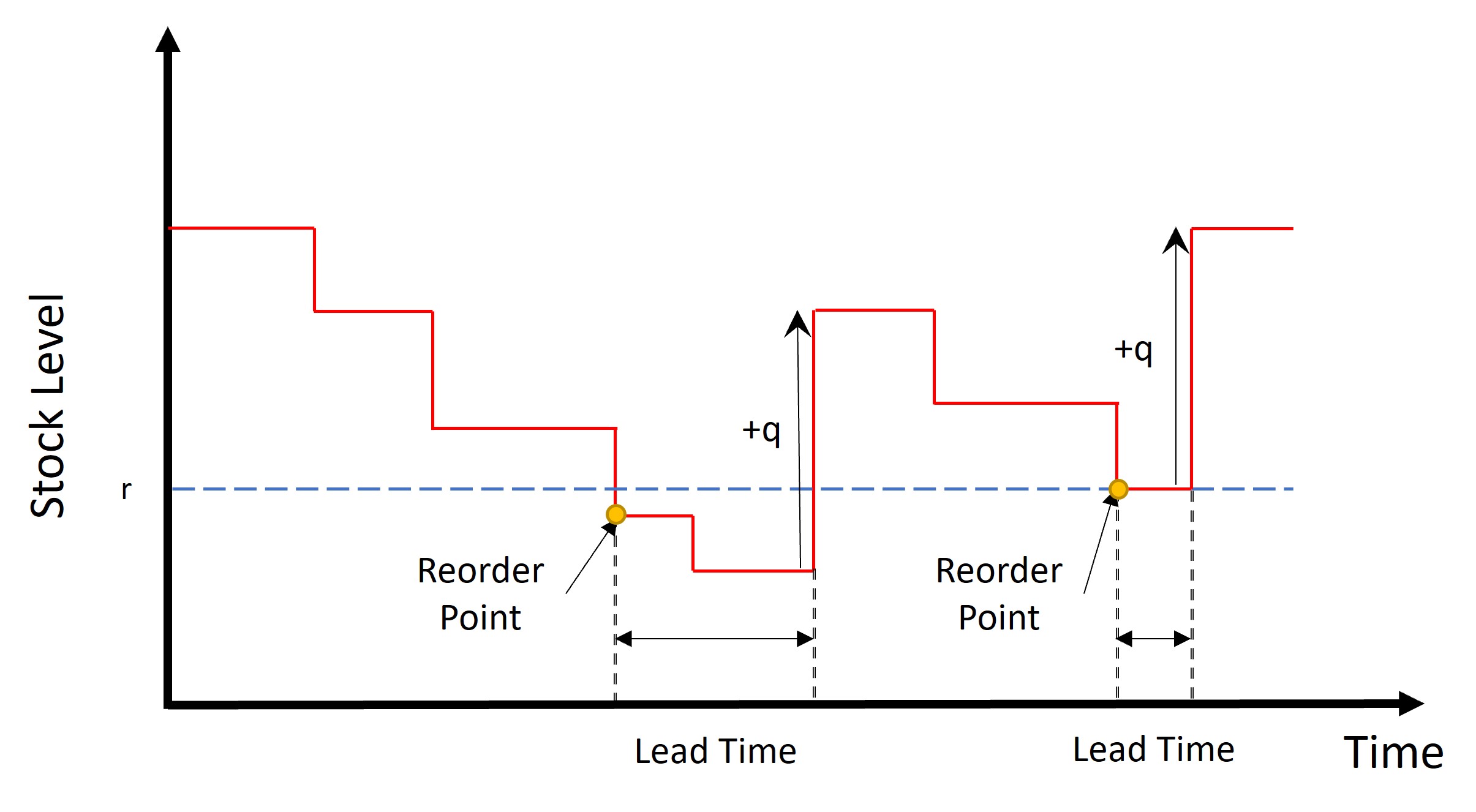}
    \caption{Stock level profile using $(r,q)$ policy}
    \label{fig:qr_inventory}
\end{figure}

The $(r, q)$ policy involves placing an order for $q$ units of stock whenever the stock level drops to or below the $r$ stock level. Subsequently, a replenishment of $q$ units of stock will be made after some lead time. Figure~\ref{fig:qr_inventory} illustrates the stock level profile under this policy.

\subsubsection{$(r,q,T)$ policy} \label{sec:single_rqt}
\begin{figure}[h]
    \centering
    \includegraphics[width=0.5\textwidth]{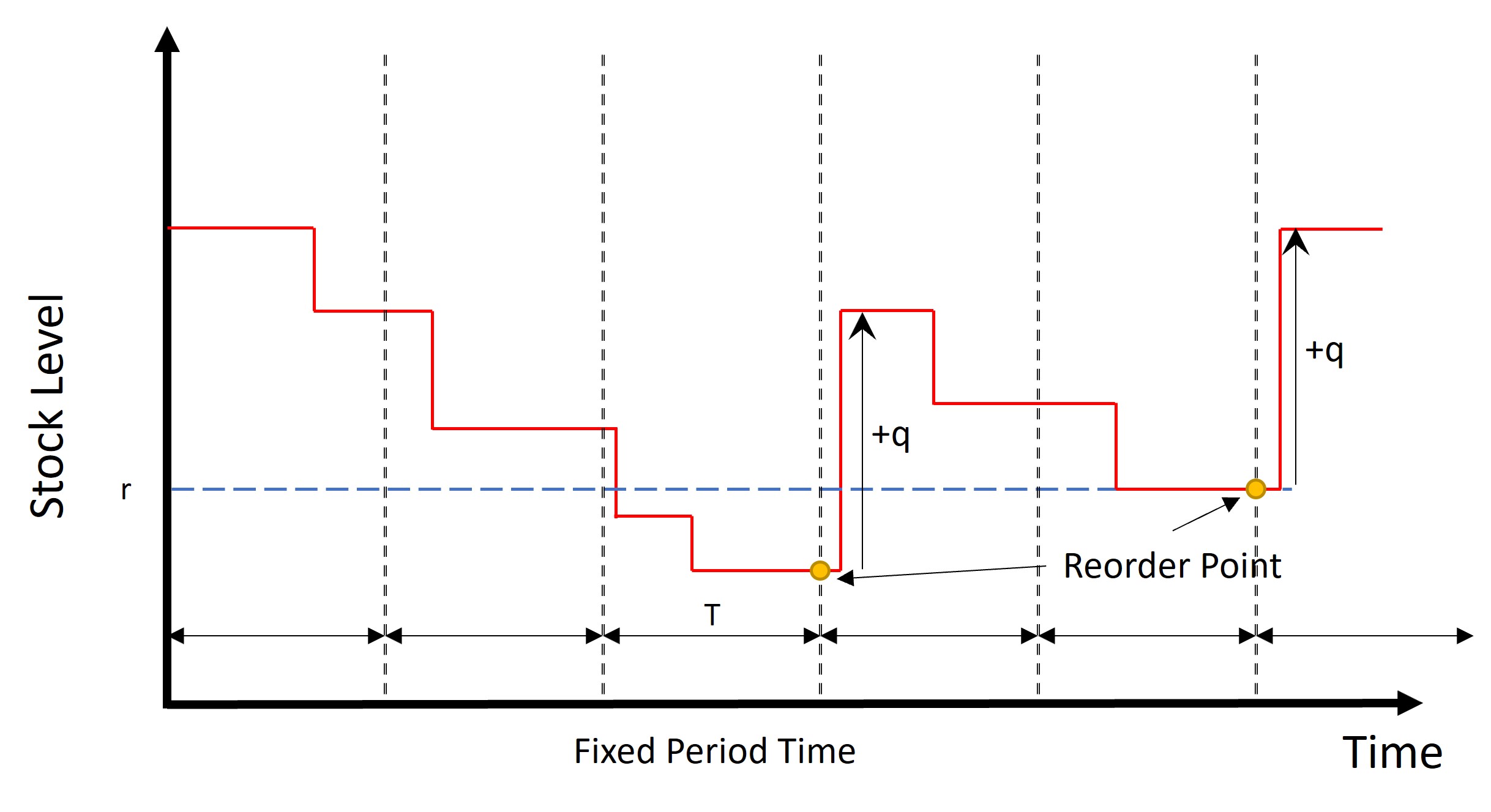}
    \caption{Stock level profile using $(r, q, T)$ policy}
    \label{fig:qrt_inventory}
\end{figure}

The $(r, q, T)$ policy involves placing an order for $q$ units of stock if the stock level drops to or below the $r$ stock level at every review period $T$. Subsequently, a replenishment of $q$ units of stock will be made after some lead time. The difference between $(r,q)$ and $(r,q,T)$ is that the former can make an order at any moment while the latter can place an order only during the review period, once the stock level is equal to or less than $r$. Figure~\ref{fig:qrt_inventory} illustrates the stock level profile under this policy.

\subsection{Markov Chain Model}
A Markov chain is a mathematical modeling method describing a sequence of state transitions with the assumption that the transition depends only on the previous state. It is used for defining the two strategies mathematically.

For example, let the random variable $X_k$ represent the number of satellites in an in-plane orbit at step $k$ with state space $\mathcal{X}= \lbrace 0, 1, \cdots, N \rbrace$. Then the probability distribution of being state levels at time step $k$ is expressed as:
\begin{equation} 
    \pi_k = \begin{bmatrix}
        P(X_k = N) & P(X_k = N-1) & \cdots & P(X_k = 0)
    \end{bmatrix}^\top
\end{equation}
and $\pi_k$ is an $N+1 \times 1$ column vector. Due to satellite failures and replenishment, there will be state transitions at each step, and their probability can be expressed as:
\begin{equation} 
    p_{ij, k} = P(X_{k+1} = i | X_{k} = j )
\end{equation}
assuming the Markov chain. With this notation, the probability distribution for the next step state can be obtained as:
\begin{equation} 
    P(X_{k+1} = i) = P(X_{k+1} = i | X_{k} = j ) P(X_{k} = j) \implies \pi_{k+1} = P_{k}\pi_{k} 
\end{equation}
where $P_{k}\in \mathbb{R}^{(N+1)\times(N+1)}$ is the state transition matrix at step $k$ and $(P_k)_{ij} = p_{ij, k}$. If the transition matrix is not a function of the step, then the distribution after the $k$ step transition can be expressed as:
\begin{equation} 
    \pi_{k} = P\pi_{k-1} = \cdots = P^{k}\pi_0
\end{equation}
where $\pi_0$ is the initial state distribution. 

\subsubsection{Stationary Distribution} The distribution is called stationary distribution for $P$ if it satisfies:
\begin{equation} 
    \pi = P \pi
\end{equation}
and it is an eigenvector of $P$ with eigenvalue $1$. If $P$ is finite, irreducible, and all its states are positive recurrent, then the Markov chain has a unique stationary distribution \cite{hillier2015introduction}. Without proof, the spare management satisfies this condition due to state decrease by failure and state increase by replenishment. On the other hand, one can interpret the meaning of the $\pi$ as the long-run proportion of time that the chain spends in states. For example, the $j$-th element of $\pi$ represents the probability that the state is $j$ at some random time in the future. 

\subsubsection{State Division Matrix} We may be interested in the portion of state distribution, e.g. the probability distribution of state level being larger than $r$ or less than or equal to $r$. For this reasons, we define the following state division matrices:
\begin{equation}
    C_{r}^+ \equiv \begin{bmatrix}
        I_{N-r} & 0 \\ 0 & 0
    \end{bmatrix},\ 
    C_{r}^- \equiv \begin{bmatrix}
        0 & 0 \\ 0 & I_{r+1}
    \end{bmatrix}
\end{equation}
which implies:
\begin{equation}
    I_{N+1} = \begin{bmatrix}
        I_{N-r} & 0 \\ 0 & I_{r+1}
    \end{bmatrix}
    = C_{r}^+ + C_{r}^-
\end{equation}
For a given general state distribution $\pi$:
\begin{equation} \label{eq:general_pi}
    \pi = \begin{bmatrix}
        \pi_N & \cdots & \pi_{r+1} & \pi_{r} & \cdots & \pi_0
    \end{bmatrix}^\top
\end{equation}
Multiplying the state division matrices gives the state distribution with a stock level greater than $r$ and less than or equal to $r$, as follows:
\begin{equation} 
\begin{aligned}
    C_{r}^+\pi &= \begin{bmatrix} \pi_N &  \cdots & \pi_{r+1} & 0 & \cdots & 0 \end{bmatrix}^\top\\ 
    C_{r}^-\pi &= \begin{bmatrix} 0 & \cdots & 0 & \pi_{r} & \cdots & \pi_0 \end{bmatrix}^\top
\end{aligned}
\end{equation}
which are keys to deriving the mathematical formula of the resupply strategy.

\subsection{Probabilistic Model}
\subsubsection{Failure Probability Distribution}
The failure of a satellite is modeled by a Poisson distribution\cite{collopy2003assigning}, i.e., an exponential distribution for the time between failures, and the spare satellites are assumed not to fail. Note that the more accurate model for the failure distribution is known to be a Weibull distribution\cite{castet2009satellite}, but the underlying assumption of immediate failure replacement makes the use of a Poisson distribution reasonable\cite{jakob2019optimal}.

Let $T_\text{mc}$ be the time step of a Markov process and $\lambda_\text{sat}$ be the failure rate of an operating satellite per $T_\text{mc}$. Then the probability of having $k$ failures when there are $n$ satellites for $T_\text{mc}$ becomes:
\begin{equation}\label{mpois_fail}
    \nu_{k,n} = 
    P(F=k| X=n) = \begin{cases}
        \frac{(n\lambda_\text{sat})^k}{k!}e^{-n\lambda_\text{sat}} & \mbox{if } n \leq N_\text{sat}\\ 
        \frac{(N_\text{sat}\lambda_\text{sat})^k}{k!}e^{-N_\text{sat}\lambda_\text{sat}} & \mbox{if } n > N_\text{sat} \\
    \end{cases}
\end{equation}
where $F$ represents the random variable for the number of failures, $0 \leq k \leq N_\text{sat}$, and $0 \leq n$. Additionally, let $\bar N_\text{sat} > N_\text{sat}$ be the maximum number of in-plane satellites, including spares. Then the state transition matrix due to failure can be defined as:
\begin{equation} \label{eq:fail_matrix_sim}
    P_f = \begin{bmatrix}
        \nu_{0, \bar N_\text{sat}}    & \cdots & 0 & \cdots & 0\\
        \vdots    & \ddots & \vdots & \ddots & \vdots \\
        1 -\sum_{k=0}^{N_\text{sat}}\nu_{k,\bar N_\text{sat}}  & \cdots & \nu_{0, N_\text{sat}}  & \cdots & 0\\
        0 & \cdots & \nu_{1,  N_\text{sat}} & \cdots & 0\\
        \vdots   & \ddots & \vdots &\ddots & \vdots \\
        0    & \cdots & 1 -\sum_{k=0}^{N_\text{sat}} \nu_{k, N_\text{sat}} & \cdots & 1 \\
    \end{bmatrix}
\end{equation}
where $P_f \in \mathbb{R}^{(\bar N_\text{sat}+1 )\times (\bar N_\text{sat}+1)}$.
Note that the column sum of the state transition matrix must be equal to $1$ for it to be proper, which is indeed the case here.

\subsubsection{Lead Time Probability Distribution} The lead time of ground resupply insertion is modeled by a shifted exponential distribution\cite{jakob2019optimal} as:
\begin{equation}
    T \sim \text{Exp}(\mu_\text{LV}) + T_\text{LV}
\end{equation}
and the probability density function of an exponential distribution is defined as:
\begin{equation}
    f(T=t;\mu) = 
    \begin{cases}
        \mu e^{-\mu t} & \quad t \geq 0\\
        0 & \quad t < 0\\
    \end{cases}
\end{equation}
where $\mu_{\text{LV}}$ is the inverse of the mean of exponential lead time distribution and $T_\text{LV}$ is constant lead time. Then the probability of having a lead time between $k$ and $k+1$ time steps of $T_\text{mc}$ is computed as
\begin{equation}
\begin{aligned}
    \rho_k &= P(kT_\text{mc} \leq T < (k+1) T_\text{mc}) \\
    &= \int_{kT_\text{mc}}^{(k+1)T_\text{mc}} f(t-T_\text{LV};\mu_\text{LV})\ dt  \\
    &= e^{-\mu_\text{LV} kT_\text{mc}}\left( 1 - e^{-\mu_\text{LV} T_\text{mc}}\right),\ kT_\text{mc} \geq T_\text{LV}
\end{aligned}
\end{equation}
Note that to avoid additional complexity, we assumed that $T_\text{mc}$ is designed such that $T_\text{LV}$ becomes an integer multiple of $T_\text{mc}$.

\section{Direct Resupply Modeling}
This section presents the analysis method of the in-plane $(r,q)$ policy introduced in the preliminaries for the direct resupply strategy. First, we introduce the following state distributions. Note that the goal is to compute $\pi^\text{dr}$.
\begin{itemize}
    \item $\pi^q$: State distribution right after $q$ replenishment.
    \item $\pi^r$: State distribution at the reorder moment.
    \item $\pi^\text{np}$: State distribution during the non-reordering period. (i.e., reorder is not made)
    \item $\pi^\text{wp}$: State distribution during the reorder waiting period. (i.e., reorder is made)
    \item $\pi^\text{dr}$: State distribution during the entire period under the strategy.
\end{itemize}
The maximum number of satellites is $\bar N_\text{sat} = q + r$, so $\pi^{(\cdot)} \in \mathbb{R}^{\bar N_\text{sat}+1}$ from the definition of $(r,q)$.

\begin{figure}[h]
    \centering
    \includegraphics[width=0.5\textwidth]{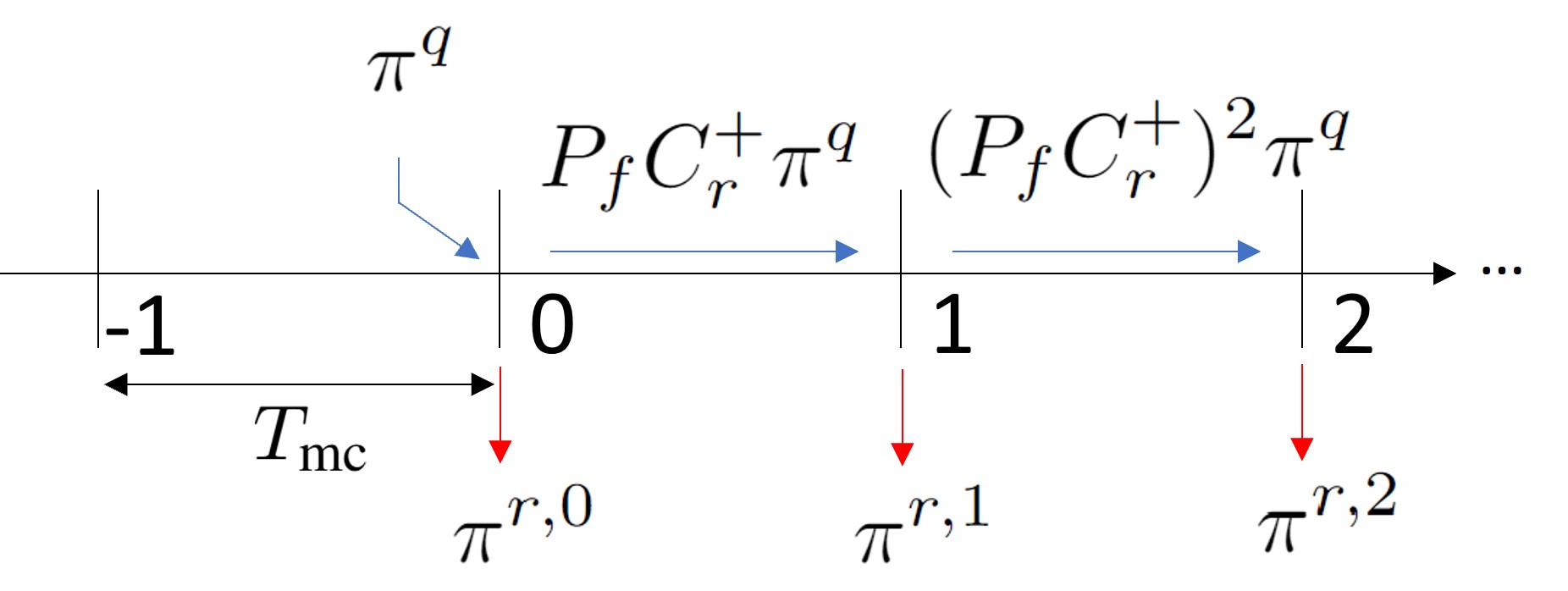}
    \caption{Propagation of state distribution during non-reordering period}
    \label{fig:direct_q2r}
\end{figure}

Let's assume the reorder has just arrived before time step $0$ as illustrated in Figure~\ref{fig:direct_q2r}. At time step $0$, the state distribution is $\pi^q$ by definition. Then, $C_{r}^+\pi^q$ will not trigger a reorder, while $C_{r}^-\pi^q$ will trigger a reorder at this time step. Likewise, the portion that has not triggered a reorder, $C_{r}^+\pi^q$, will undergo failure, so the state distribution becomes $P_f C_{r}^+ \pi^q$ at time step $1$. Then, the portion of the distribution that will trigger a reorder is $C_{r}^- P_f C_{r}^+ \pi^q$, while the portion that will not trigger a reorder is $C_{r}^+ P_f C_{r}^+ \pi^q$ at time step $1$.

Let $\pi^{r,k}$ be the relative state distribution at time step $k$ given that a reorder is made at time step $0$. From this observation, it is derived that:
\begin{equation}
\begin{aligned}
    \pi^{r,0} &= C_r^- \pi^q\\
    \pi^{r,1} &= C_r^- P_f C_r^+ \pi^q \\
    &\vdots \\
    \pi^{r,k} &= C_r^{-}(P_fC_r^+)^k \pi^q\\
    &\vdots \\
\end{aligned}
\end{equation}
One can see that each event is mutually exclusive (i.e., reordering cannot occur at different time steps simultaneously), and if we consider all cases as $k\to\infty$, the union of them is collectively exhaustive. Also, the $L_1$ norm of the relative state distribution $\pi^{r,k}$ represents the probability of the event, implying that the expected state distribution at the reorder point $\pi^r$ is computed as:
\begin{equation} \label{eq:direct_q2r}
    \pi^{r} = \sum_{k=0}^\infty \pi^{r,k} = C_r^-(I + P_fC_r^+ + \cdots)\pi^q = C_r^-\left(I - P_f C_r^+ \right)^{-1}\pi^q = P_{r/q}\pi^q
\end{equation}
where $P_{r/q}$ is the transition matrix from $\pi^q$ to $\pi^r$. Note that the infinite sum has the analytic expression since the matrix norm of $P_fC_r^+$ is less than one.

On the other hand, suppose the event of ground resupply insertion has just occurred, and the stock level was $X = r$ before the resupply. Then, the stock level right after replenishment is $X = r+q$. Similarly, if the stock level was $X=r+1$ right before the ground resupply insertion, then the stock level right after the replenishment is also $X=r+1$ because such an event cannot happen, thus not changing the stock level. Based on this observation, the replenishment transition matrix $P_q$ is defined  as follows:
\begin{equation} \label{eq:resupply_trans_matrix}
    P_q = \begin{bmatrix}
        \begin{bmatrix} I_q \\ 0_{l \times q} \end{bmatrix} & I_{r+1} \\
        0_{q \times q} & 0_{q \times (r+1)} 
        \end{bmatrix}\ \text{if}\ q \leq r+1,
    \quad 
    P_q = \begin{bmatrix}
            I_q & \begin{bmatrix} I_{r+1} \\ 0_{l\times (r+1)}  \end{bmatrix} \\
            0_{(r+1)\times q}   & 0_{(r+1)\times (r+1)}
        \end{bmatrix}\ \text{otherwise}
\end{equation}
where $l = \lvert q - (r+1) \rvert$ and $P_q \in \mathbb{R}^{(\bar N_\text{sat}+1) \times (\bar N_\text{sat}+1)}$.

\begin{figure}[h]
    \centering
    \includegraphics[width=0.5\textwidth]{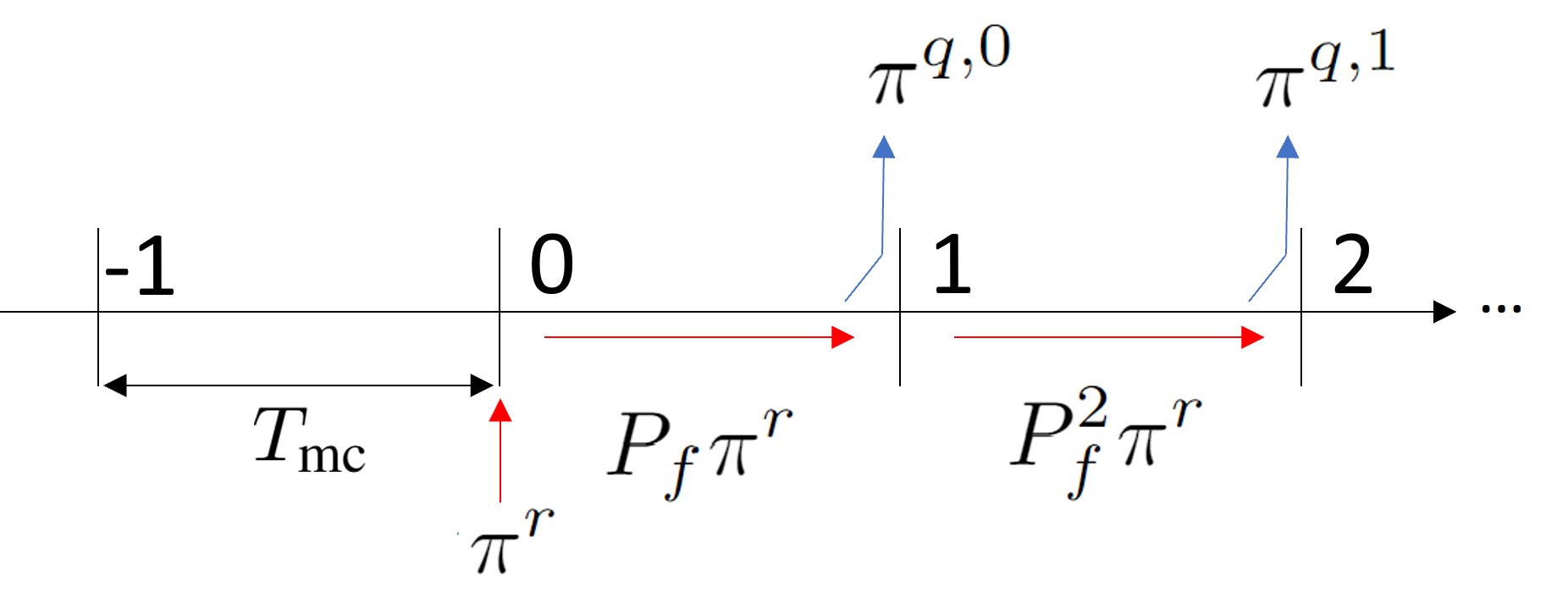}
    \caption{Propagation of state distribution during reorder waiting period}
    \label{fig:direct_r2q}
\end{figure}

Next, we will explain the transition relationship from $\pi^r$ to $\pi^q$. First, we assume that the replenishment will occur right before the time step. Consider the scenario where the reorder is made at time step $0$ with the expected state distribution $\pi^r$, as illustrated in Figure~\ref{fig:direct_r2q}. With a given lead time distribution, the probability of resupply arriving before time step $1$ is $\rho_0$ by definition. If that happens, the portion of the state distribution after receiving the spares between time steps $0$ and $1$ becomes $\rho_0 P_q P_f \pi^r$. Similarly, if the spares are received between time steps $1$ and $2$, it is $\rho_1 P_q P_f^2 \pi^r$.

Let $\pi^{q,k}$ be the relative state distribution which had a replenishment between time steps $k$ and $k+1$. Then, based on the previous investigation, it is computed as follows:
\begin{equation}
\begin{aligned}
    \pi^{q,0} &= \rho_0 P_q P_f \pi^r \\
    \pi^{q,1} &= \rho_1 P_q P_f^2 \pi^r \\
    &\vdots \\
    \pi^{q,k} &= \rho_k P_q P_f^{k+1} \pi^r\\
    &\vdots \\
\end{aligned}
\end{equation}
From the same reasoning made above, each event is mutually exclusive and collectively exhaustive, and the expected state distribution right after the replenishment $\pi^q$ is computed as:
\begin{equation} \label{eq:direct_r2q}
    \pi^{q} = \sum_{k=0}^\infty \pi^{q,k} = P_q P_f(\rho_0 I + \rho_1 P_f + \cdots)\pi^r = P_{q/r}\pi^r
\end{equation}
where $P_{q/r}$ is the transition matrix from $\pi^r$ to $\pi^q$. Note that the truncated summation should generally be applied, but there is an analytic expression for some distributions, and the exponential distribution is one of them. In the case of the exponential distribution with constant bias, the equation becomes:
\begin{equation}
    \pi^{q} = P_{q/r}\pi^r = (1-e^{-\mu_\text{LV} T_\text{mc}})P_qP_f^{m+1}(I-e^{-\mu_\text{LV} T_\text{mc}}P_f )^{-1}\pi^r
\end{equation}
where $m$ is $m = T_\text{LV}/T_\text{mc}$. 

From the state transition matrix obtained in Eqs. \eqref{eq:direct_q2r} and \eqref{eq:direct_r2q}, the self-transition matrix is computed as follows:
\begin{equation}
\pi^q = P_{q/r}\pi^r,\quad \pi^r = P_{r/q}\pi^q \implies \pi^q = P_{q/r}P_{r/q}\pi^q,\quad \pi^r = P_{r/q}P_{q/r}\pi^r
\end{equation}
This self-transition matrix satisfies the necessary condition for the uniqueness of the stationary distribution\footnote{To be precise, the truncated self-transition matrix satisfies the condition}. Therefore, applying any numerical algorithm will find the stationary distribution efficiently. Recall that the $\pi^q$ and $\pi^r$ are not the distribution at each time step but it is a distribution when some event occurs (i.e. the replenishment and reorder). We need one last step to derive the average stationary distribution for the entire period, which is the one needed for the policy analysis.

Let's assume $\pi^q$ is computed from the previous equation. Then, the relative state distribution at each time step after replenishment during the non-reordering period is derived as follows:
\begin{itemize}
    \item Distribution at time step $0$ given that reorder is not triggered: $\pi^{\text{np},0} = C_r^+\pi^q$
    \item Distribution at time step $1$ given that reorder is not triggered: $\pi^{\text{np},1} = C_r^+P_fC_r^+\pi^q$
    \item Distribution at time step $k$ given that reorder is not triggered: $\pi^{\text{np},k} = C_r^+(P_fC_r^+)^k\pi^q$
\end{itemize}
Therefore, the expected state distribution during the entire non-reordering period becomes:
\begin{equation} \label{eq:direct_pi_np}
\pi^\text{np} = \frac{1}{T^\text{np}}\sum_{k=0}^\infty \pi^{\text{np},k} = \frac{1}{T^\text{np}} C_r^+(I+P_fC_r^+\cdots)\pi^q =
\frac{1}{T^\text{np}} C_r^+(I - P_fC_r^+)^{-1}\pi^q
\end{equation}
where $T^\text{np}$ is the average duration of the non-reordering period in units of $T_\text{mc}$ and is the sum of the $L_1$ norm of the relative state distribution, or equivalently summation of all elements, as:
\begin{equation}
T^\text{np} = \sum_{k=0}^\infty \lVert \pi^{\text{np},k} \rVert_{1} = \sum_{k=0}^\infty \sum_{j=0}^{\bar N_\text{sat}}
 \pi^{\text{np},k}_j
\end{equation}

Similarly, given the $\pi^r$ distribution, the relative state distribution at each time step after the trigger of reorder during the reorder waiting period is derived as follows:
\begin{itemize}
    \item Distribution at time step $0$ given that replenishment is not made: $\pi^{\text{wp},0} = \pi^r$
    \item Distribution at time step $1$ given that replenishment is not made: $\pi^{\text{wp},1} = (1-\rho_0)P_f\pi^r$
    \item Distribution at time step $k$ given that replenishment is not made: $\pi^{\text{wp},k} = (1-\sum_{i=0}^{k-1}\rho_i)P_f^k\pi^r$
\end{itemize}
Therefore, the expected state distribution during the entire reorder waiting period becomes:
\begin{equation}
    \pi^\text{wp} = \frac{1}{T^\text{wp}}\sum_{k=0}^\infty \pi^{\text{wp},k} = \frac{1}{T^\text{wp}} (I+\rho_0^c P_f + \rho_1^c P_f^2+\cdots)\pi^r
\end{equation}
where $\rho_k^c = (1-\sum_{i=0}^{k}\rho_i)$ and $T^\text{wp}$ is the average duration of the reorder waiting period in units of $T_\text{mc}$. Similarly, $T^\text{wp}$ is the sum of the $L_1$ norm of the relative state distribution:
\begin{equation}
    T^\text{wp} = \sum_{k=0}^\infty  \lVert \pi^{\text{wp},k}\rVert_{1}
\end{equation}
Based on the assumption of the lead time model, the infinite sum becomes the following explicit equation:
\begin{equation}
    \pi^\text{wp} = \frac{1}{T^\text{wp}} \left( \sum_{i=0}^{m} P_f^i + \rho_0^c P_f^{m+1}(I - e^{-\mu_\text{LV} T_\text{mc}}P_f)^{-1} \right) \pi^r
\end{equation}

Lastly, since we know the state distribution during the non-reordering period and reorder waiting period, and the length of each period, the average state distribution for the entire period of the direct resupply policy is given by:
\begin{equation}
\pi^\text{dr} = \frac{T^\text{np}}{T^\text{np} + T^\text{wp}}\pi^\text{np} + \frac{T^\text{wp}}{T^\text{np} + T^\text{wp}}\pi^\text{wp}
\end{equation}
and the average period for one replenishment cycle is:
\begin{equation}
    T_\text{cycle}^\text{dr} = \left(T^\text{np} + T^\text{wp}\right) T_\text{mc} 
\end{equation}

With this result, one can compute any performance metric as needed. For example, the average probability of having a stock level less than $y$ is given by:
\begin{equation}
    P(X < y) = \sum_{k=0}^{y-1} \pi^\text{dr}_k
\end{equation}
where $\pi_k$ is the probability having $k$ stock level as defined in Eq.\eqref{eq:general_pi}.

\section{Indirect Resupply Modeling}
The idea of analyzing the indirect resupply strategy is to solve the in-plane and parking orbits independently and then check the consistency of the results. Specifically, the in-plane stock level depends on the availability of the parking orbit, and the stock level will determine the demand for in-plane spare satellites. They are coupled through the in-plane spares demand and the parking spares availability. We used the fixed-point iteration to find a consistent solution, as will be discussed in the subsequent subsection. 

Lastly, the notation used in this section is identical to that of the direct resupply strategy except for the subscripts $(\cdot)_i$ and $(\cdot)_p$ indicating the in-plane orbit and parking orbit, respectively. In addition, the maximum number of in-plane satellites is $\bar N_{\text{sat}_i} = q_i + r_i$ in units of satellites, and the maximum number of parking spares is $\bar N_{\text{sat}_p} = q_p + r_p$ in units of batches. That is, $\pi^{(\cdot)_i} \in \mathbb{R}^{\bar{N}_{\text{sat}_i}+1}$ and $\pi^{(\cdot)_p} \in \mathbb{R}^{\bar N_{\text{sat}_p}+1}$. 
\begin{itemize}
    \item $\pi^{q_{(\ast)}}$: State distribution of in-plane/parking orbit right after $q_{(\ast)}$ replenishment.
    \item $\pi^{r_{(\ast)}}$: State distribution of in-plane/parking orbitat the reorder moment.
    \item $\pi^\text{np$_{(\ast)}$}$: State distribution of in-plane/parking orbit during the non-reordering period.
    \item $\pi^\text{wp$_{(\ast)}$}$: State distribution of in-plane/parking orbit during the reorder waiting period. 
    \item $\pi^\text{ir$_{(\ast)}$}$: State distribution of in-plane/parking orbit during the entire period under the strategy.
\end{itemize}

\subsection{In-Plane Orbit Markov Chain}
Based on the assumptions made in the preliminary section, the time period for an in-plane orbit to align with the subsequent parking orbit can be computed as follows:
\begin{equation} \label{eq:t_plane_indirect}
    T_\text{plane} = \frac{2\pi}{N_\text{park}|\dot \Omega_\text{plane} - \dot \Omega_\text{park}|}
\end{equation}
By properly choosing the $T_\text{mc}$ value, we can make $T_\text{plane}$ an integer multiple of $T_\text{mc}$ as
\begin{equation}
    T_\text{plane} = k_i T_\text{mc},\quad k_i\in \mathbb{N}
\end{equation}
and $T_\text{plane}$ becomes the review period of in-plane $(r,q,T)$ policy. Note that the failure transition matrix for the unit time step is identical to that of direct resupply. Therefore, if the state distribution right after the replenishment (i.e., RAAN contact) was $\pi^{q_i}$ then the state distribution $\pi^{r_i}$ becomes
\begin{equation} \label{eq:indirect_qi2ri}
    \pi^{r_i} = P_{f_i} \pi^{q_i} = P_{f}^{k_i} \pi^{q_i}
\end{equation}

Let's assume the stock level of each parking orbit follows independent and identically distributed conditions. Also, assume that the expected state distribution of the parking orbit right before the contact point is known. Then we can define the following parking availability probability, which will be computed in the parking analysis part, as:
\begin{equation} \label{eq:parking_avail}
    \kappa_j = P(X_p \geq D_i = j)=\sum_{k=j}^{\bar N_{\text{sat}_p}} \pi^{\text{ir}_p}_k, \quad j\in\mathbb{N},\quad \kappa_0 = 1
\end{equation}
where $\pi^{\text{ir}_p}_k = P^{\text{ir}_p}(X_p = k)$ is an element of $\pi^{\text{ir}_p}$,
$\bar N_{\text{sat}_p}$ is the maximum number of parking stock, and $D_i$ is the random variable of demand for spare satellites in the in-plane at the contact moment, which is defined as:
\begin{equation} \label{eq:park_demand_func}
    D_{i} = \begin{cases}
        \lceil\frac{ r_i + 1 - X_i}{q_i}\rceil & \text{if}\ X_i \leq r_i \\
        0 & \text{if}\ X_i > r_i
    \end{cases}
\end{equation}
where $\lceil \cdot \rceil$ is the ceiling operator. With this definition, the in-plane replenishment transition matrix becomes:
\begin{equation} \label{eq:park_avail_mat}
    P_{q_i} = \begin{bmatrix}
        \kappa_0I_{q_i} & \kappa_1 I_{q_i} & \kappa_2 I_{q_i}  & \cdots \\
        0 & (\kappa_0 -\kappa_1)I_{q_i} & (\kappa_1 -\kappa_2) I_{q_i}  & \cdots \\
        0 & 0 & (\kappa_0 -\kappa_1) I_{q_i}  & \cdots \\
        \vdots & \vdots & \vdots  & \ddots
    \end{bmatrix}
\end{equation}
Note that the Eq.~\eqref{eq:park_demand_func} and \eqref{eq:park_avail_mat} have this expression because they are measured in units of batch $q_i$.

From the state transition matrices $P_{f_i}$ and $P_{q_i}$, the self-transition matrices are computed as follows:
\begin{equation} \label{eq:inplane_sol_eq}
\pi^{q_i} = P_{q_i}P_{f_i} \pi^{q_i},\quad \pi^{r_i} = P_{f_i}P_{q_i}\pi^{r_i}
\end{equation}
These self-transition matrices satisfy the necessary condition for the existence of a unique stationary distribution. Therefore, any numerical algorithm applied will efficiently find the stationary distribution. With this stationary distribution, we can compute the average stock level during the replenishment cycle as
\begin{equation}
    \pi^{\text{ir}_i} = \frac{1}{k_i}\left( I + P_f + P_f^2 + \cdots + P_f^{k_i-1} \right)\pi^{q_i}
\end{equation}
In addition, one can compute the in-plane demand probability distribution using $\pi^{r_i}$ as:
\begin{equation} \label{eq:inplane_demand_prob}
\eta_{d} = P(D_i = d) = \sum_{k=dq_i}^{(d+1)q_i-1} \pi^{r_i}_{\bar N_{\text{sat}_i} - k}, \quad d\in \mathbb{N}
\end{equation}
where $\pi^{r_i}_{j} = P^{r_i}(X_i = j)$ is an element of $\pi^{r_i}$ and the parking orbit will receive this demand probability distribution at every RAAN contact moment (i.e., every $T_\text{park}$).

\subsection{Parking Orbit Markov Chain}
Now, we will explain the analysis method for the parking orbits. First, the time period for a parking orbit to align with the subsequent in-plane orbit can be computed as follows:
\begin{equation} \label{eq:t_park_indirect}
    T_\text{park} = \frac{2\pi}{N_\text{plane}|\dot \Omega_\text{plane} - \dot \Omega_\text{park}|} = k_p T_\text{mc},\quad k_p \in \mathbb{N}
\end{equation}
Unlike the case of the in-plane orbit, the stock level of parking orbits only drops at the contact moment by sending the spares to the in-plane orbit. Therefore, the stock transition matrix of the parking orbit between contact moments is:
\begin{equation} \label{eq:fail_matrix_park}
    P_{f_p} = \begin{bmatrix}
        \eta_{0} & 0   & \cdots & 0\\
        \eta_{1} & \eta_{0}   & \cdots & 0 \\
        \vdots & \vdots   & \ddots & \vdots \\
        1 -\sum_{i=0}^{N_p}\eta_{i} & 1 -\sum_{i=0}^{N_p-1}\eta_{i}   & \cdots & 1 \\
    \end{bmatrix}
\end{equation}
In addition, we assume that the request for ground resupply only happens at the contact moment, so that $T_\text{park}$ becomes the review period of parking orbit $(r,q,T)$ policy. Similar to the case of the direct resupply strategy, the replenishment transition matrix, $P_{q_p}\in \mathbb{R}^{(\bar N_{\text{sat}_p} + 1) \times (\bar N_{\text{sat}_p} + 1)}$,  is defined as Eq.~\eqref{eq:resupply_trans_matrix} with $q_p$ and $r_p$.

Lastly, the demand from the in-plane orbit occurs only every $T_\text{park}$ period, whereas the ground resupply can arrive at any moment within the $T_\text{park}$ period (with a time step of $T_\text{mc}$). In addition, the reorder of parking spares will only be triggered at the contact moment, and the replenishment will not be made until $T_\text{LV}$ has passed since the reorder trigger moment. This complicated structure requires additional steps for deriving the formula.

\begin{figure}[h]
    \centering
    \includegraphics[width=0.65\textwidth]{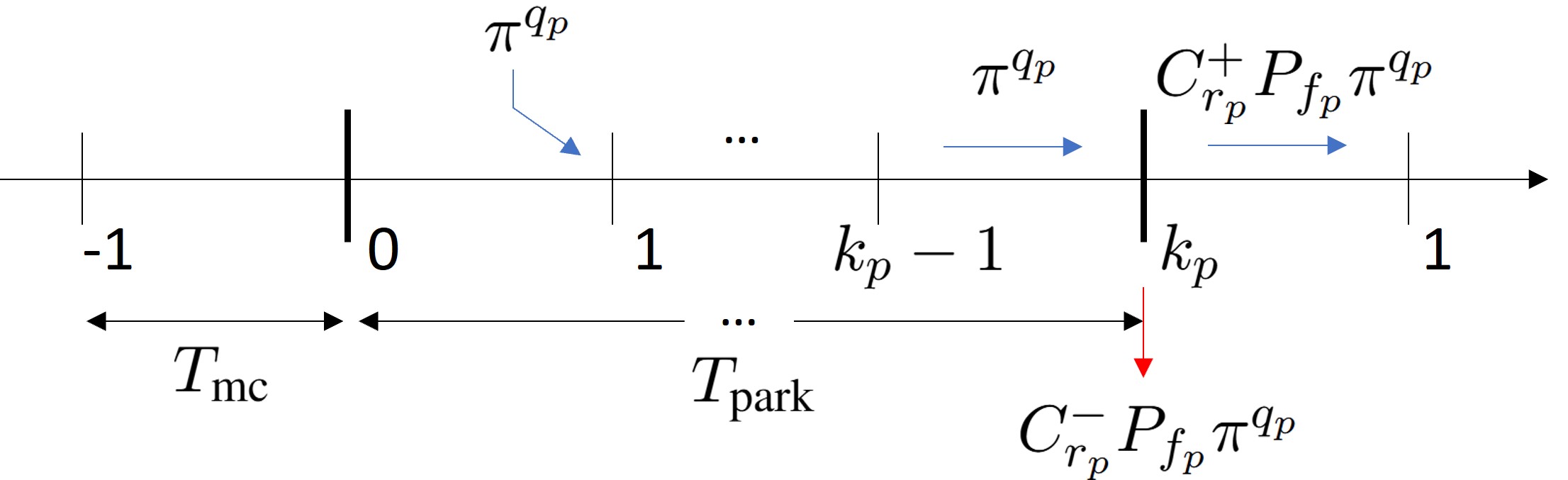}
    \caption{Propagation of state distribution during non-reordering period}
    \label{fig:indirect_park_q2r}
\end{figure}

Let's assume the reorder has arrived at a time step within $0$ to $k_p$ and the state distribution at the step is $\pi^{q_p}$, as shown in Figure~\ref{fig:indirect_park_q2r}. Then, $C_{r_p}^+ P_{f_p}\pi^{q_p}$ will not trigger a reorder, while $C_{r_p}^- P_{f_p}\pi^{q_p}$ will trigger a reorder at $0$ step (RAAN match at 0 or $k_p$ time step). Note that the structure is identical to the state transition of the direct resupply, hence the following relationship holds
\begin{equation} \label{eq:indirect_q2r_park}
    \pi^{r_p} =C_{r_p}^-(P_{f_p} + P_{f_p}^2 C_{r_p}^+ + \cdots)\pi^{q_p} = C_{r_p}^-P_{f_p}\left(I - P_{f_p} C_{r_p}^+ \right)^{-1}\pi^{q_p} 
\end{equation}

\begin{figure}[h]
    \centering
    \includegraphics[width=0.65\textwidth]{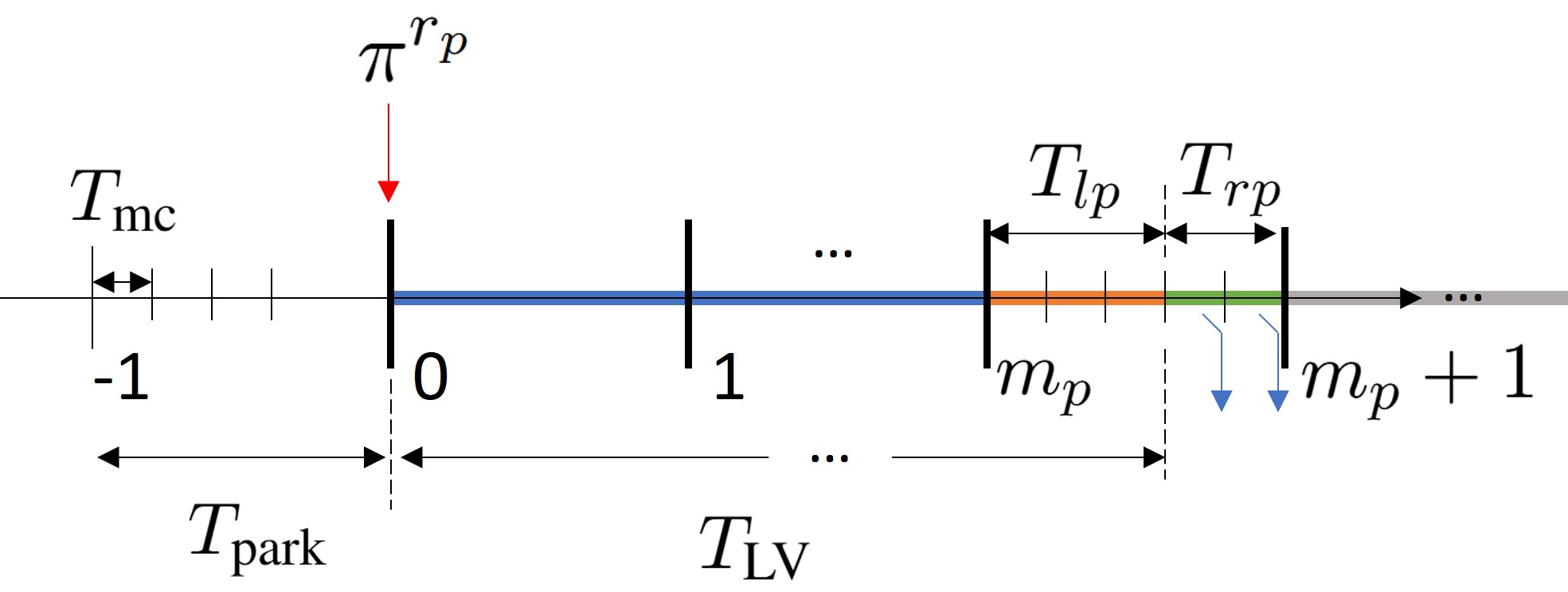}
    \caption{Propagation of state distribution during waiting period}
    \label{fig:indirect_park_r2q}
\end{figure}

Likewise, let's assume the reorder has been made at step 0 and the state distribution at the step is $\pi^{r_p}$. For the shifted exponential lead time distribution, we have to consider four terms (colored segments): 
the distribution during the quotient and remainder parts when $T_\text{LV}$ is divided by $T_\text{park}$, and the tail end part. Referring to Figure~\ref{fig:indirect_park_r2q}, the relevant parameters are defined as:
\begin{equation}
    m_p = \left\lfloor \frac{T_\text{LV}}{T_\text{park}} \right\rfloor, \quad 
    T_{lp} = T_\text{LV} - m_p T_\text{park} ,\quad
    T_{rp} = (m_p+1)T_\text{park} - T_\text{LV}
\end{equation}
From the beginning of the third part(green), there is a probability of making a replenishment. The state distribution at the beginning is $P_{f_p}^{m_p}\pi^{r_p}$ since it only undergoes the failure (i.e., demand for in-plane), and the cumulative probability of the third part and the fourth part are:
\begin{equation}
    \rho_{p_3} = 1 - e^{\mu_\text{LV} T_{rp}},\quad \rho_{p_4} = e^{\mu_\text{LV} T_{rp}} (1 - e^{-\mu_\text{LV} T_\text{park}}) 
\end{equation}
Therefore, the state distribution right after replenishment becomes:
\begin{equation} \label{eq:indirect_r2q_park}
    \pi^{q_p} = P_{q_p} P_{f_p}^{m_p} \left( \rho_{p_3}I + \rho_{p_4}P_{f_p} 
    \left(I - e^{-\mu_\text{LV} T_\text{park}} P_{f_p} \right)^{-1} \right)\pi^{r_p}
\end{equation}
Rearranging the Eqs.~\eqref{eq:indirect_q2r_park} and \eqref{eq:indirect_r2q_park} gives the self-transition relationship, so both $\pi^{q_p}$ and $\pi^{r_p}$ can be obtained.

Next is to compute the average distribution during the waiting period $\pi^{\text{wp}_p}$ and the non-reordering period $\pi^{\text{np}_p}$. First, $\pi^{\text{wp}_p}$ can be computed as a weighted sum of the distribution of the four segments as:
\begin{equation}
    \pi^{\text{wp}_p} = \frac{1}{T^{\text{wp}_p}}\sum_{k=1}^{4} \pi^{\text{wp}_p,k}, 
    \quad T^{\text{wp}_p} = \big\lVert \sum_{k=1}^4 \pi^{\text{wp}_p,k} \big\rVert_{1}
\end{equation}
where
\begin{equation}
\begin{aligned}
    \pi^{\text{wp}_p,1} &= k_p \left(I+P_{f_p} + \cdots + P_{f_p}^{m_p-1} \right)\pi^{r_p}\\
    \pi^{\text{wp}_p,2} &= \frac{T_{lp}}{T_\text{mc}} P_{f_p}^{m_p}\pi^{r_p}\\
    \pi^{\text{wp}_p,3} &= \frac{ e^{-\mu_\text{LV} T_\text{mc}} \cdot ( e^{-\mu_\text{LV} T_{rp}} - 1) }{ e^{-\mu_\text{LV} T_\text{mc}} - 1}P_{f_p}^{m_p}\pi^{r_p}\\
    \pi^{\text{wp}_p,4} &= \frac{ e^{-\mu_\text{LV} (T_{rp} + T_\text{mc})} \cdot 
    ( e^{-\mu_\text{LV} T_\text{park}} - 1) }{ e^{-\mu_\text{LV} T_\text{mc}} - 1}
    P_{f_p}^{m_p+1}\left(I - e^{-\mu_\text{LV} T_\text{park}} P_{f_p} \right)^{-1}    
    \pi^{r_p}\\
\end{aligned}
\end{equation}
There is no resupply during the first $m_p$ cycles, so the stock distribution will keep decreasing during this period, as expressed in the first term. The stock level is constant during the second period, and the normalized duration is $T_{lp}$ over $T_\text{mc}$. During the third and fourth segments, the portion of the state will experience replenishment, and this can be expressed as the sum of a geometric series. Note that the equation for the sum of a finite geometric series is used in the third and fourth terms:
\begin{equation}
    a + ar + ar^2 + \cdots + ar^{n-1} = \frac{a\cdot (r^n - 1)}{r-1}
\end{equation}

The average state distribution for the non-reordering period is almost identical to the Eq.~\eqref{eq:direct_pi_np} except for the additional term due to arbitrary arrival time of ground resupply within $T_\text{park}$ period, which is derived as:
\begin{equation} \label{eq:indirect_pi_np_park}
    \pi^{\text{np}_p} = \frac{1}{T^{\text{np}_p}} 
    \left(
    k_p C_{r_p}^+P_{f_p}(I - C_{r_p}^+P_{f_p})^{-1}
    + \frac{1}{T_\text{mc}}\sum_{l=0}^{k_p-1} \bar T_l \cdot \bar \rho_l
    \right)\pi^q 
\end{equation}
where the average non-reordering period is computed as:
\begin{equation}
    T^{\text{np}_p} = 
    \big\lVert
    k_p C_{r_p}^+P_{f_p}(I - C_{r_p}^+P_{f_p})^{-1}
    + \frac{1}{T_\text{mc}}\sum_{l=0}^{k_p-1} \bar T_l \cdot \bar \rho_l
    \pi^q  \big \rVert_1
\end{equation}
, and the additional variables are defined as
\begin{equation}
    \bar T_l = \begin{cases}
        T_{rp} - lT_\text{mc} & \mbox{if } T_{rp} - lT_\text{mc} \geq 0 \\
        T_\text{park} + T_{rp} - lT_\text{mc} & \mbox{if } T_{rp} - lT_\text{mc} < 0 \\
    \end{cases}
\end{equation}
and
\begin{equation}
    \bar \rho_l = \frac{1 - e^{-\mu_\text{LV} T_\text{mc}}}{1-e^{-\mu_\text{LV} T_\text{park}}} e^{-l \mu_\text{LV} T_\text{mc}}
\end{equation}
Referring to Figure~\ref{fig:indirect_park_q2r}, if replenishment is made right before time step 1, then the state distribution will not change up to time step $k_p-1$ (this is equivalent to the $T_{rp}$ period). On the other hand, if it happens right before time step 0, then the state distribution will change at that time step. In addition, the probability of having replenishment at $k$ time steps from the beginning of $T_{rp}$ is known from the exponential distribution. Then the expectation can be expressed as the second term of Eq.~\eqref{eq:indirect_pi_np_park}.

Lastly, since we know the state distribution during the non-reordering period and reorder waiting period, and the length of each period, the average state distribution for the entire period of the parking orbit of indirect resupply policy is given by:
\begin{equation} \label{eq:indirect_park_sol}
\pi^{\text{ir}_p} = \frac{T^{\text{np}_p}}{T^{\text{np}_p} + T^{\text{wp}_p}}\pi^{\text{np}_p} + \frac{T^{\text{wp}_p}}{T^{\text{np}_p} + T^{\text{wp}_p}}\pi^{\text{wp}_p}
\end{equation}
and the average period for one replenishment cycle is:
\begin{equation}
    T_\text{cycle}^{\text{ir}_p} = \left(T^{\text{np}_p} + T^{\text{wp}_p}\right) T_\text{mc} 
\end{equation}
and this is used for updating Eq.~\eqref{eq:parking_avail}.

\subsection{Overall Procedure to Compute the Solution}
We divide the analysis of coupled in-plane and parking orbits into two parts, each explained in the previous subsections. Solving each analysis independently will yield inconsistent results, necessitating iterative adjustments to achieve a consistent solution. The adjusting process uses the parking availability probability from the parking orbit analysis Eq.~\eqref{eq:parking_avail} and the in-plane spare demand probability from the in-plane analysis Eq~\eqref{eq:inplane_demand_prob}. The detailed procedure for solving the indirect resupply method is summarized in Table.

\begin{algorithm} 
    \caption{Fixed Point Iteration for the Indirect Resupply}
  \begin{algorithmic}[1]
    \REQUIRE Constellation Configuration: $(a_i,i,N_\text{sat}, N_\text{plane})$, Failure Probability Model: $\nu_{ij}$,
    LV Lead Time Model: $\rho_i'$, Minimum Time Step: $T_\text{sim}$
    \INPUT Parking Orbit Parameter: $(a_p, N_\text{park})$, In-Plane Resupply Parameter: $(r_i, q_i)$, Parking Resupply Parameter: $(r_p, q_p)$, Time Step: $T_\text{mc}$.
    \OUTPUT Average In-plane Distribution $\pi^{\text{ir}_i}$, Average Parking Distribution $\pi^{\text{ir}_p}$
    \STATE Compute $T_\text{plane}$, $T_\text{park}$ using Eq.~\eqref{eq:t_plane_indirect} and Eq.~\eqref{eq:t_park_indirect}
    \STATE Compute $P_{f_i}$ Eq.~\eqref{eq:indirect_qi2ri} and Eq.~\eqref{eq:fail_matrix_sim}.
    \STATE Set $\kappa = [1\ 1\ \cdots\ 1]^\top$ and Initialize $P_{q_i}$ as Eq.~\eqref{eq:park_avail_mat}
    \WHILE{$ |\eta^{k+1} - \eta^{k}| > \varepsilon $ or $ |\kappa^{k+1} - \kappa^{k}| >\varepsilon $ and $k < k^{\max}$}
      \STATE Compute in-plane orbit solution $\pi^{q_i}$ and $\pi^{r_i}$ using Eq.~\eqref{eq:inplane_sol_eq}
      \STATE Update in-plane orbit demand probability $\eta$ using Eq.~\eqref{eq:inplane_demand_prob}
      \STATE Compute the parking orbit demand transition matrix $P_{f_p}$ using Eq.~\eqref{eq:fail_matrix_park}
      \STATE Compute parking orbit solution $\pi^{q_p}$ and $\pi^{r_p}$ using Eq.~\eqref{eq:indirect_q2r_park} and Eq.~\eqref{eq:indirect_r2q_park}
      \STATE Compute average parking orbit solution $\pi^{\text{ir}_p}$ using Eq.~\eqref{eq:indirect_park_sol}
      \STATE Update parking availability probability $\kappa$ using Eq.\eqref{eq:parking_avail}
      \STATE Compute the in-plane replenishment matrix $P_{q_i}$ using Eq.\eqref{eq:park_avail_mat}
      \STATE $k = k + 1$
    \ENDWHILE
  \end{algorithmic}
\end{algorithm}
It turns out that the solution converges with a relative error of less than $10^{-5}$ within 10 iterations for the practical problem.

\section{Validation and Optimization Result}
\subsection{Validation of Direct Resupply Strategy Analysis Method}
The following parameters are used for the validation of direct resupply policy.
\begin{table}[H]
    \centering
    \begin{tabular}{c|c|c}
    Parameters & Value  & Description \\ \hline
    $q$ & 4  & The number of spares per replenishment $[\text{\#}]$ \\
    $r$ & 42 & The reorder condition of a direct replenishment $[\text{\#}]$ \\
    $N_\text{sat}$ & 40 & Nominal number of satellites per in-plane orbit $[\text{\#}]$ \\
    $\lambda_\text{sat}$ & 0.05 - 0.15 & Satellite failure rate per year per satellite $[\text{\#}]$ \\
    $T_\text{mc}$ & 1 & Markov Process (and Simulation) time step $[\text{day}]$ \\
    $1/\mu_\text{LV}$ & 60 & Average lead time of launch vehicle $[\text{day}]$ \\
    $T_\text{LV}$ & 30 & Constant lead time of launch vehicle for in-plane orbit $[\text{day}]$ \\
    \hline
    \end{tabular}
    \caption{Parameters for the direct resupply method analysis}
\end{table}
We have validated the analysis method under various failure rates. The blue bar represents the average state distribution from the Monte Carlo simulations, and the red star indicates the stationary solution from the proposed method. The tested failure rates are $0.05$, $0.1$, and $0.15$, and they are referred to as low, moderate, and high failure scenarios, respectively.

\begin{figure}[!h]
    \centering
    \includegraphics[width=.32\textwidth]{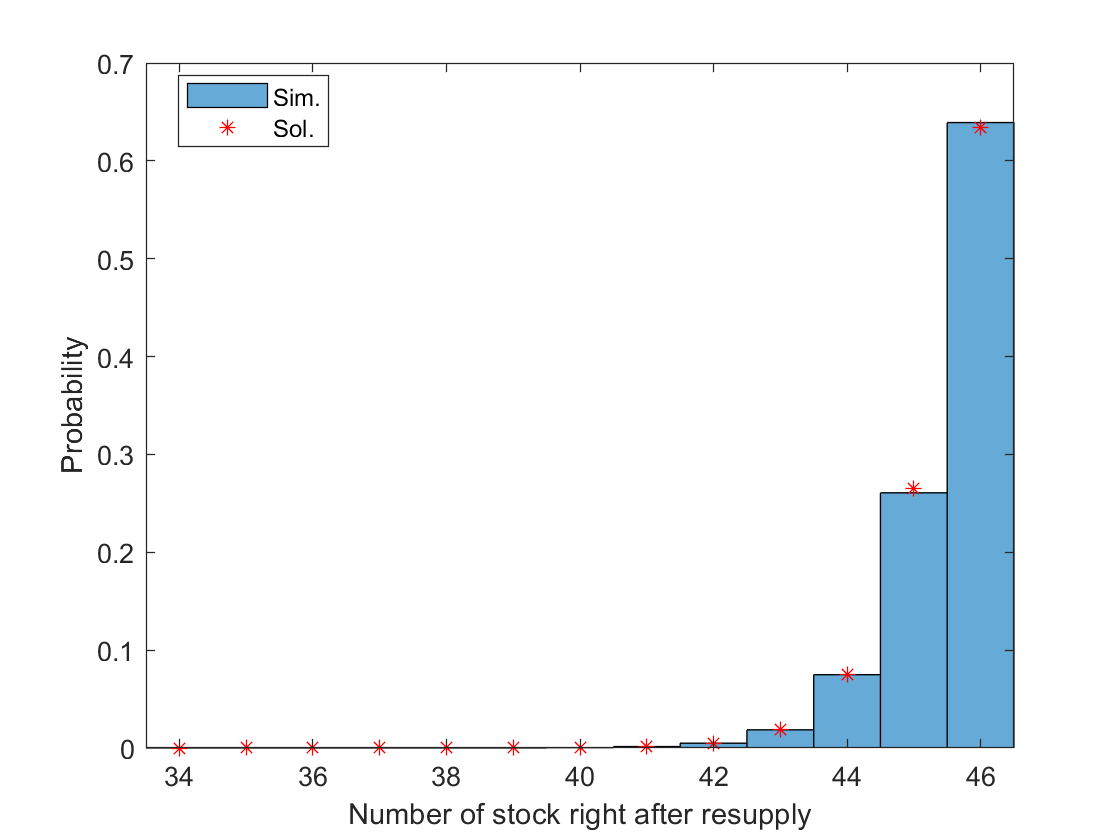}
    \includegraphics[width=.32\textwidth]{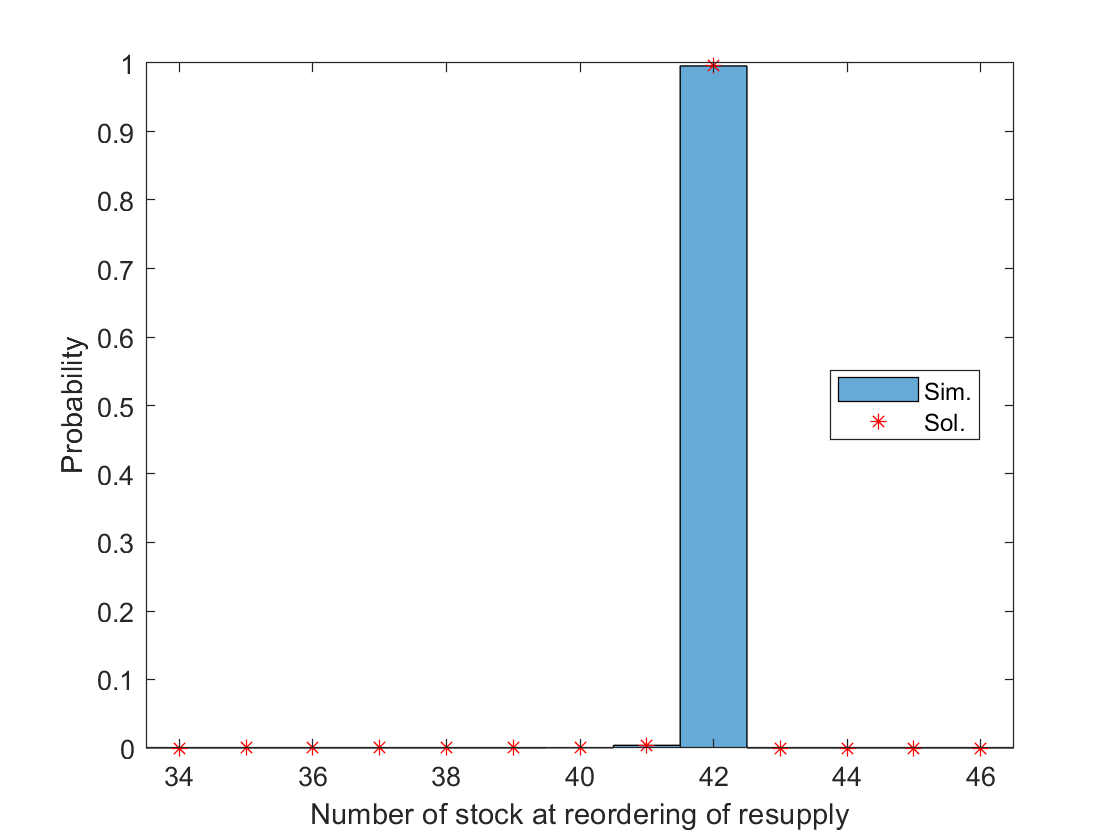}
    \includegraphics[width=.32\textwidth]{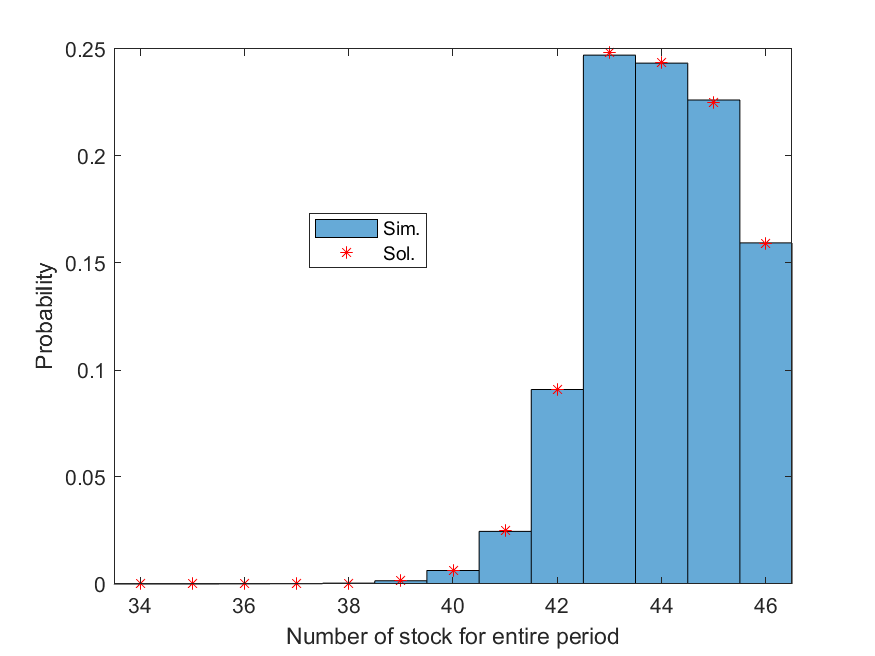}
    \caption{Stationary solution of (a) $\pi^{q}$ (b) $\pi^{r}$ (c) $\pi^{\text{dr}}$ for low failure rate}
\end{figure}

\begin{figure}[!h]
    \centering
    \includegraphics[width=.32\textwidth]{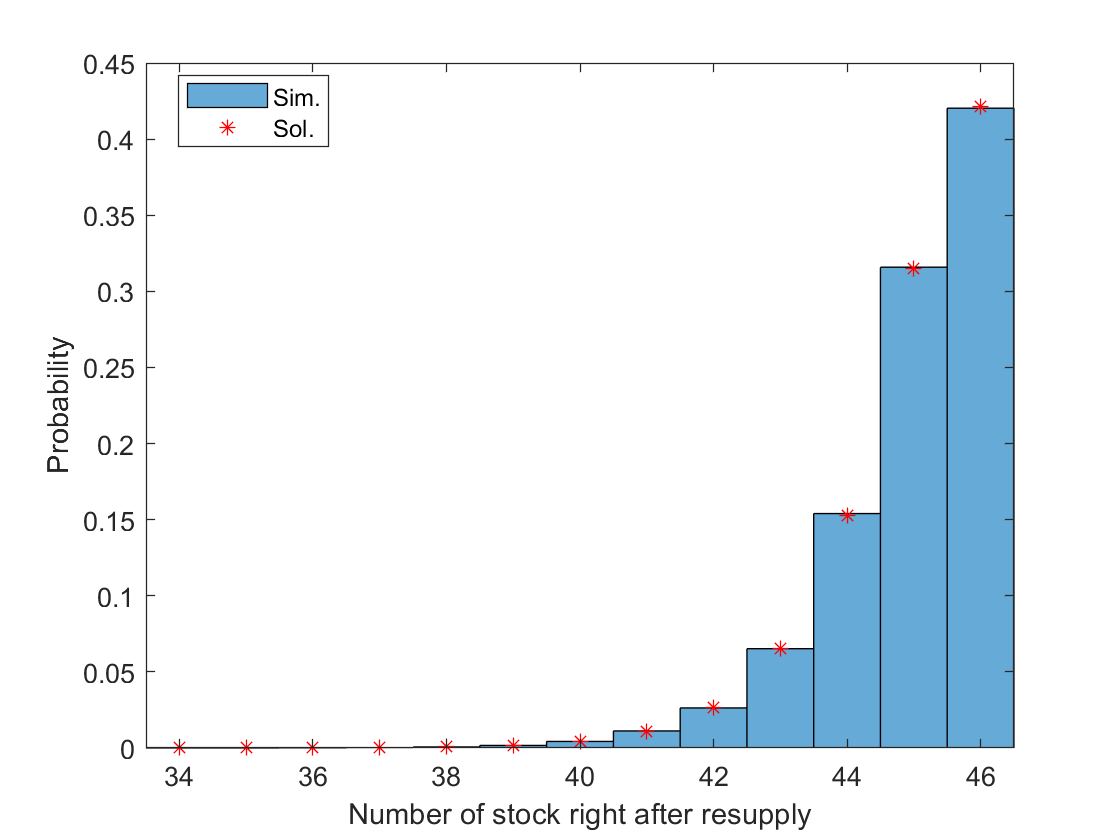}
    \includegraphics[width=.32\textwidth]{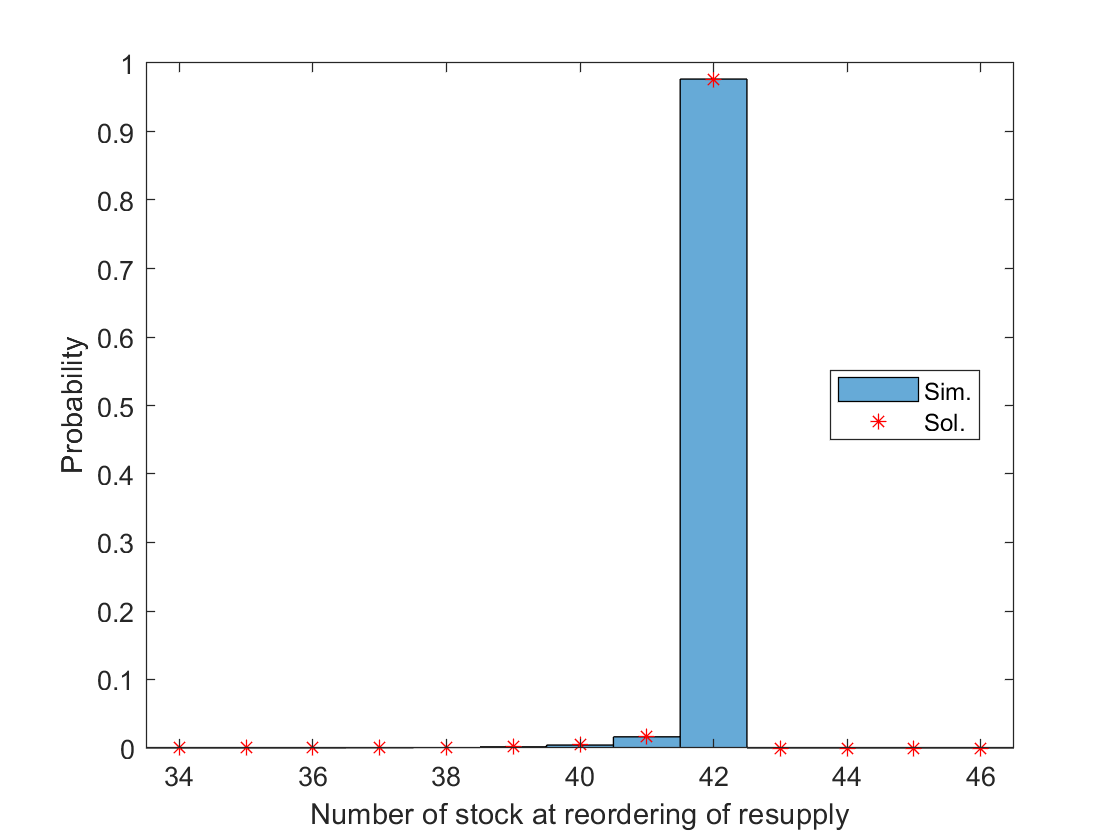}
    \includegraphics[width=.32\textwidth]{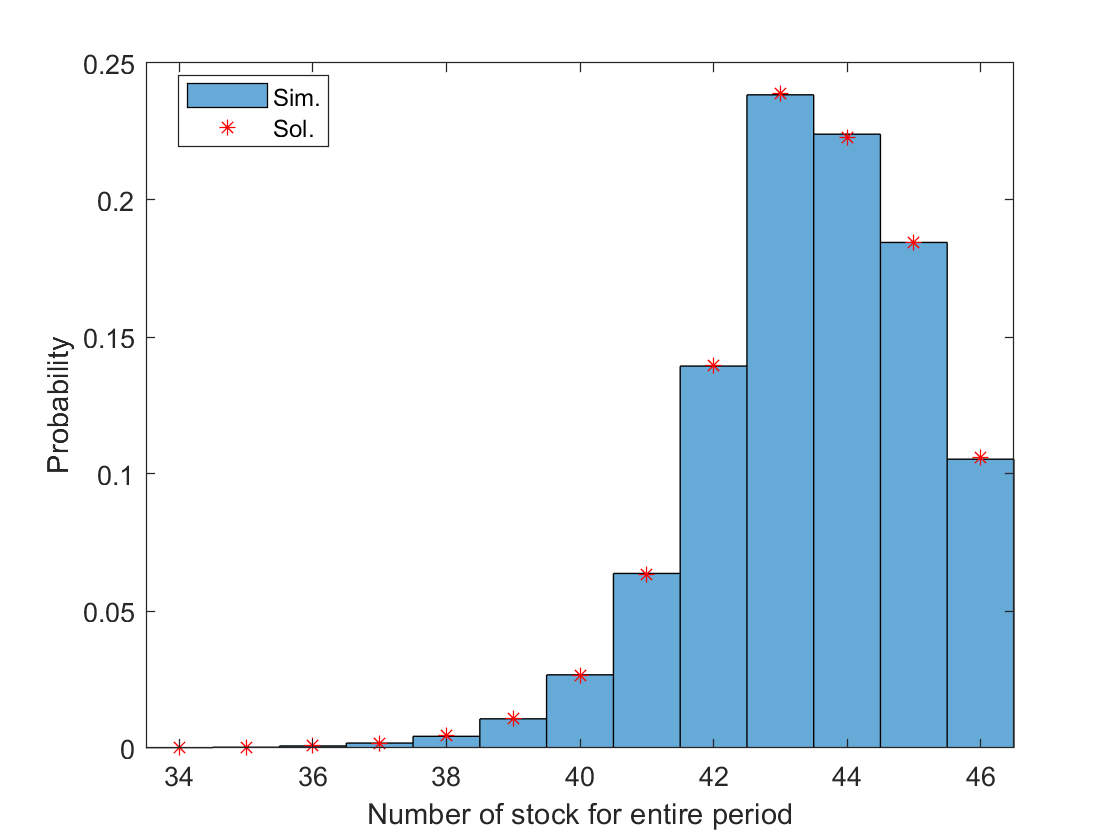}
    \caption{Stationary solution of (a) $\pi^{q}$ (b) $\pi^{r}$ (c) $\pi^{\text{dr}}$ for moderate failure rate}
\end{figure}

\begin{figure}[!h]
    \centering
    \includegraphics[width=.32\textwidth]{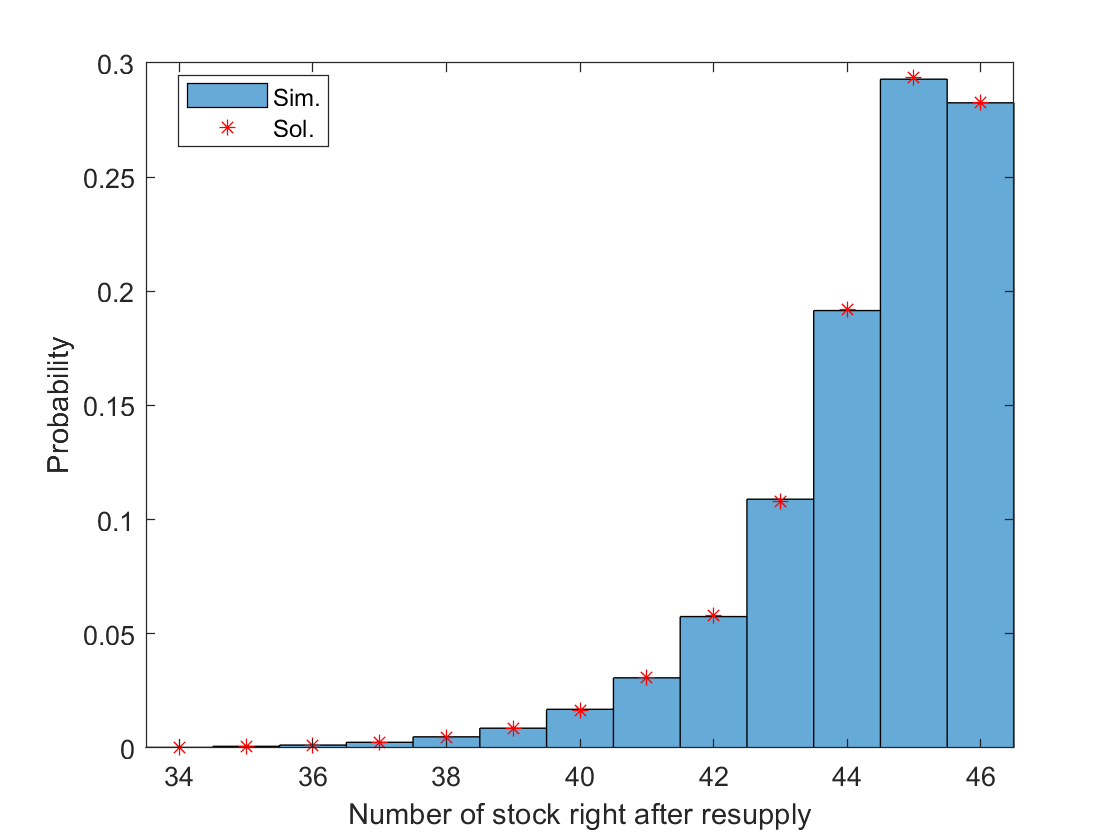}
    \includegraphics[width=.32\textwidth]{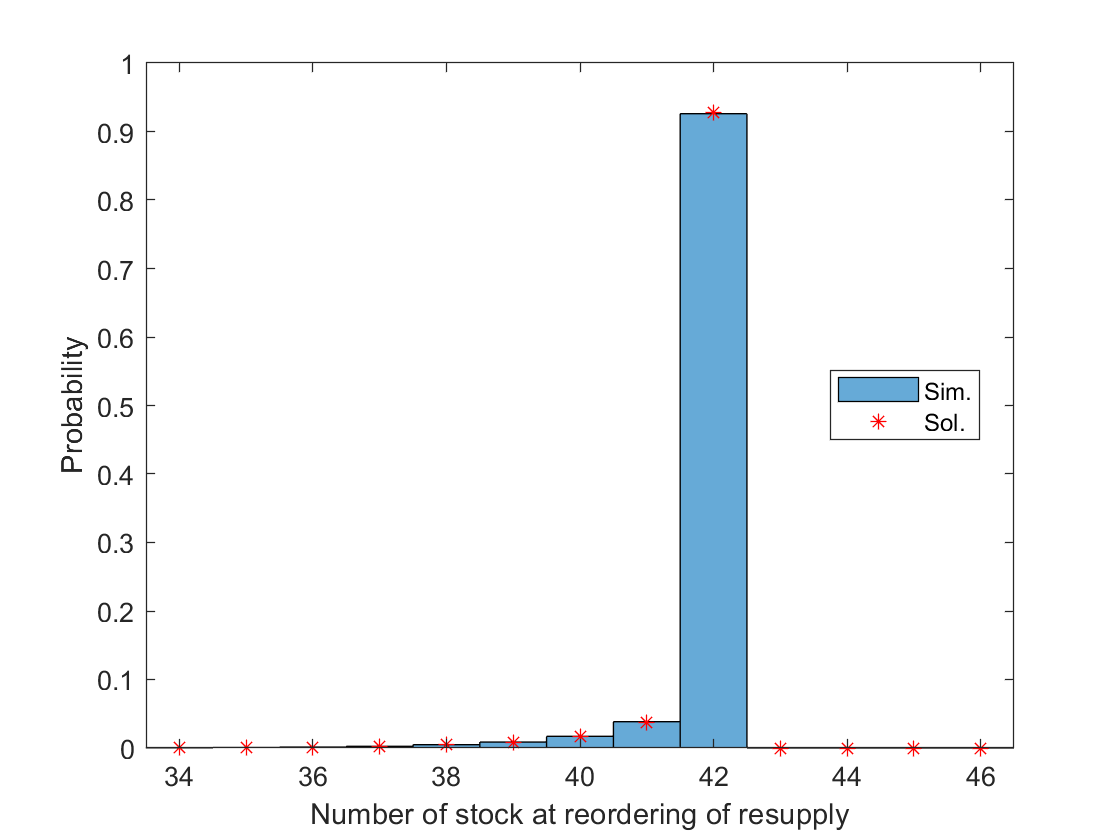}
    \includegraphics[width=.32\textwidth]{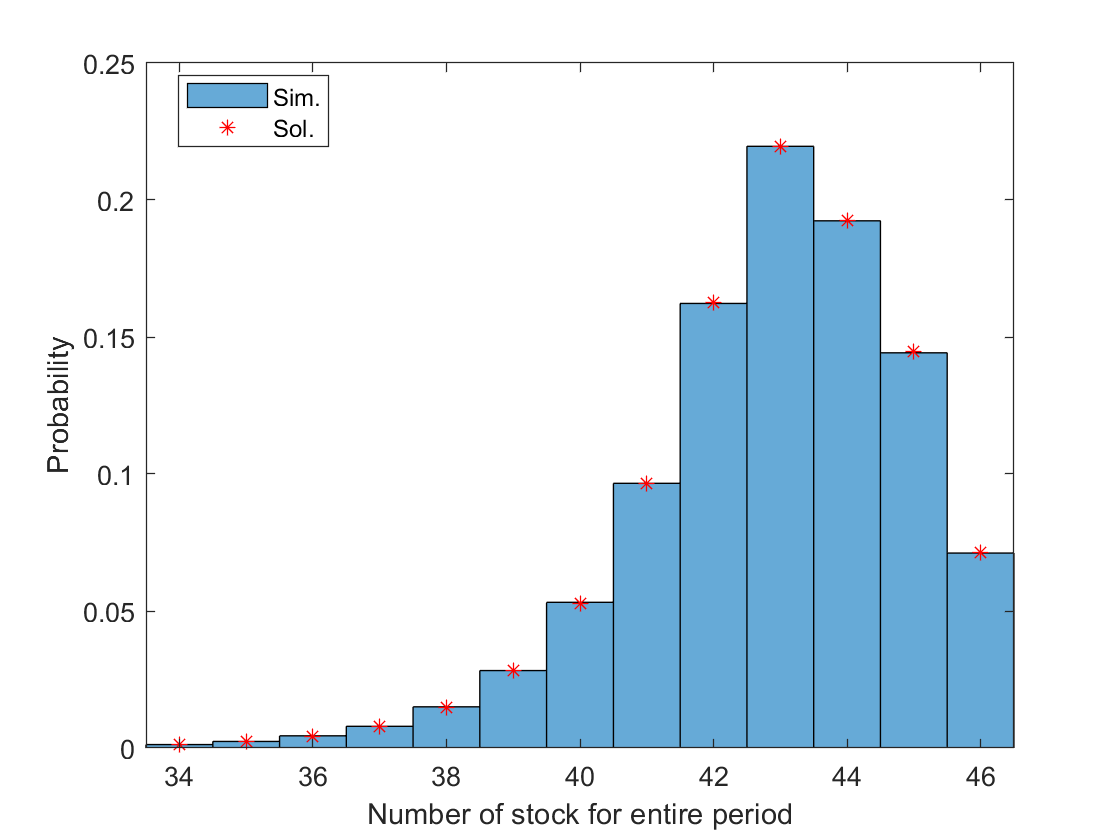}
    \caption{Stationary solution of (a) $\pi^{q}$ (b) $\pi^{r}$ (c) $\pi^{\text{dr}}$ for high failure rate}
\end{figure}
Based on the result figures, the proposed Markov process analysis accurately computes the expected state distribution. In addition, a higher failure rate while fixing the $(r,q)$ parameters results in a decrease in state level, as expected. Note that the proposed method only requires milliseconds, while the simulation takes around 30 minutes in a typical personal computing environment. Due to its fast but accurate computation, the proposed method can be integrated into the optimization problem, as will be discussed in a later section.

\subsection{Validation of Indirect Resupply Analysis Method}
The following parameters are used for the validation of indirect resupply policy.
\begin{table}[H]
    \centering
    \begin{tabular}{c|c|c}
    Parameters & Value  & Description \\ \hline
    $q_i$ & 4  & The number of spares satellites per replenishment for in-plane orbit $[\text{\#}]$ \\
    $r_i$ & 42 & The reorder condition of a replenishment for in-plane orbit $[\text{\#}]$ \\
    $q_p$ & 8  & The number of spares batches per replenishment for parking orbit $[\text{\#}]$ \\
    $r_p$ & 8 & The reorder condition of a replenishment for parking orbit $[\text{\#}]$ \\
    $N_\text{sat}$ & 40 & Nominal number of satellites per in-plane orbit $[\text{\#}]$ \\
    $\lambda_\text{sat}$ & 0.05 - 0.15 & Satellite failure rate per year per satellite $[\text{\#}]$ \\
    $N_\text{plane}$ & 40 & The number of in-plane orbits $[\text{\#}]$ \\
    $T_\text{plane}$ & 200 & RAAN match period for in-plane orbits $[\text{day}]$ \\
    $N_\text{park}$ & 3 & The number of parking orbits $[\text{\#}]$ \\
    $T_\text{park}$ & 15 & RAAN match period for parking orbits $[\text{day}]$ \\
    $T_\text{mc}$ & 1 & Markov Process (and Simulation) time step $[\text{day}]$ \\
    $1/\mu_\text{LV}$ & 60 & Average lead time of launch vehicle for in-plane orbit $[\text{day}]$ \\
    $T_\text{LV}$ & 30 & Constant lead time of launch vehicle for in-plane orbit $[\text{day}]$ \\
    \hline
    \end{tabular}
    \caption{Parameters for the indirect resupply method analysis}
\end{table}
Like previous direct resupply method validation, the analysis method for the indirect method is validated under various failure rates. 
\begin{figure}[!h]
    \centering
    \includegraphics[width=.32\textwidth]{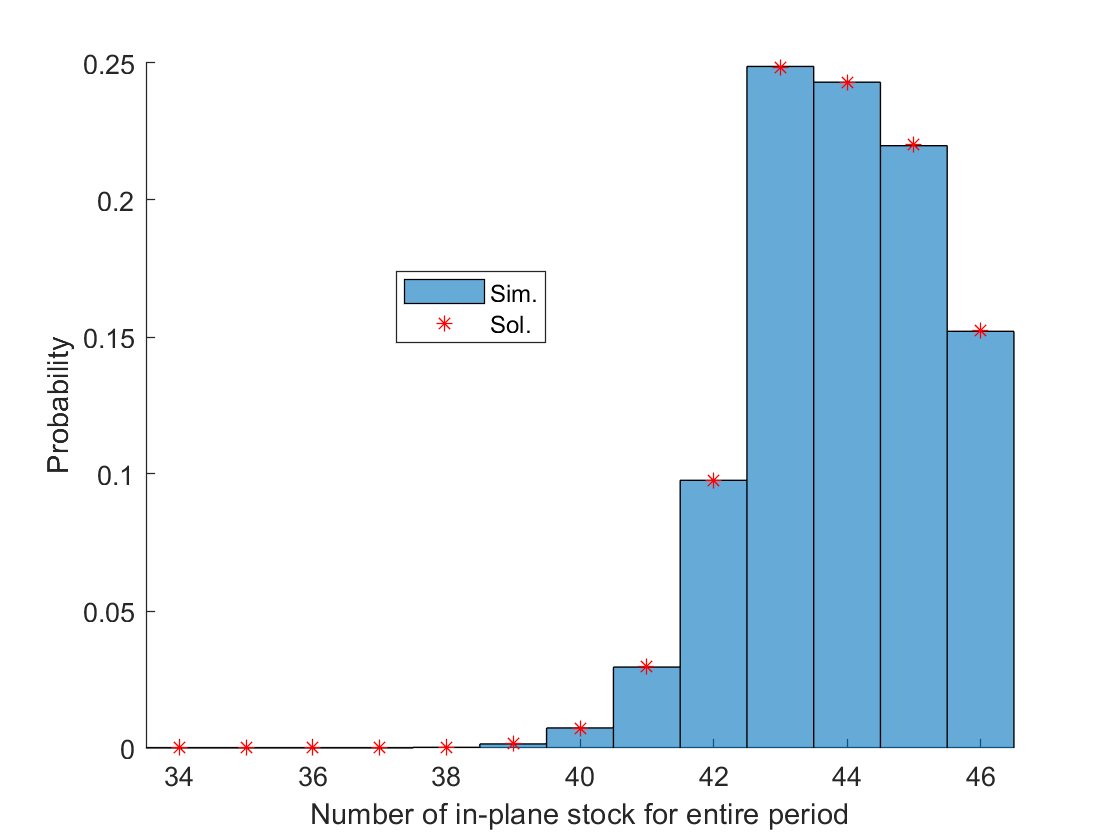}
    \includegraphics[width=.32\textwidth]{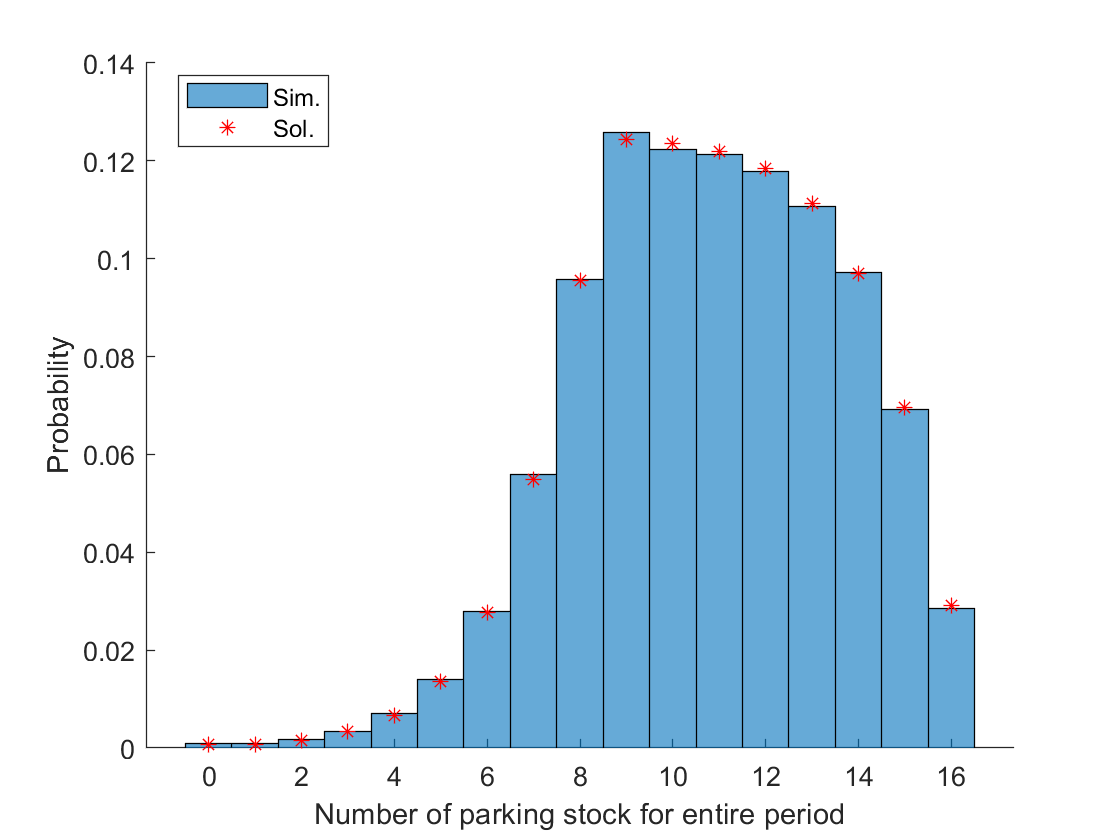}
    \includegraphics[width=.32\textwidth]{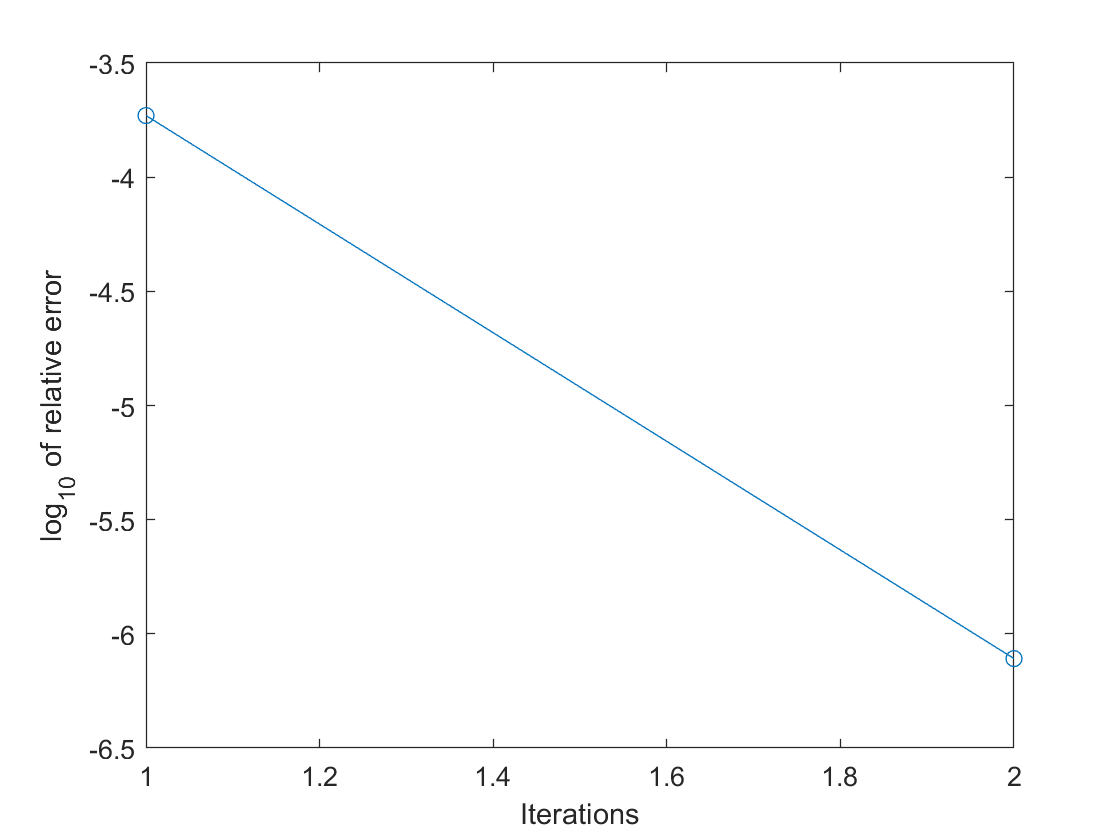}
    \caption{Stationary solution of (a) $\pi^{\text{ir}_i}$ (b) $\pi^{\text{ir}_p}$ and (c) relative error for low failure rate}
\end{figure}

\begin{figure}[!h]
    \centering
    \includegraphics[width=.32\textwidth]{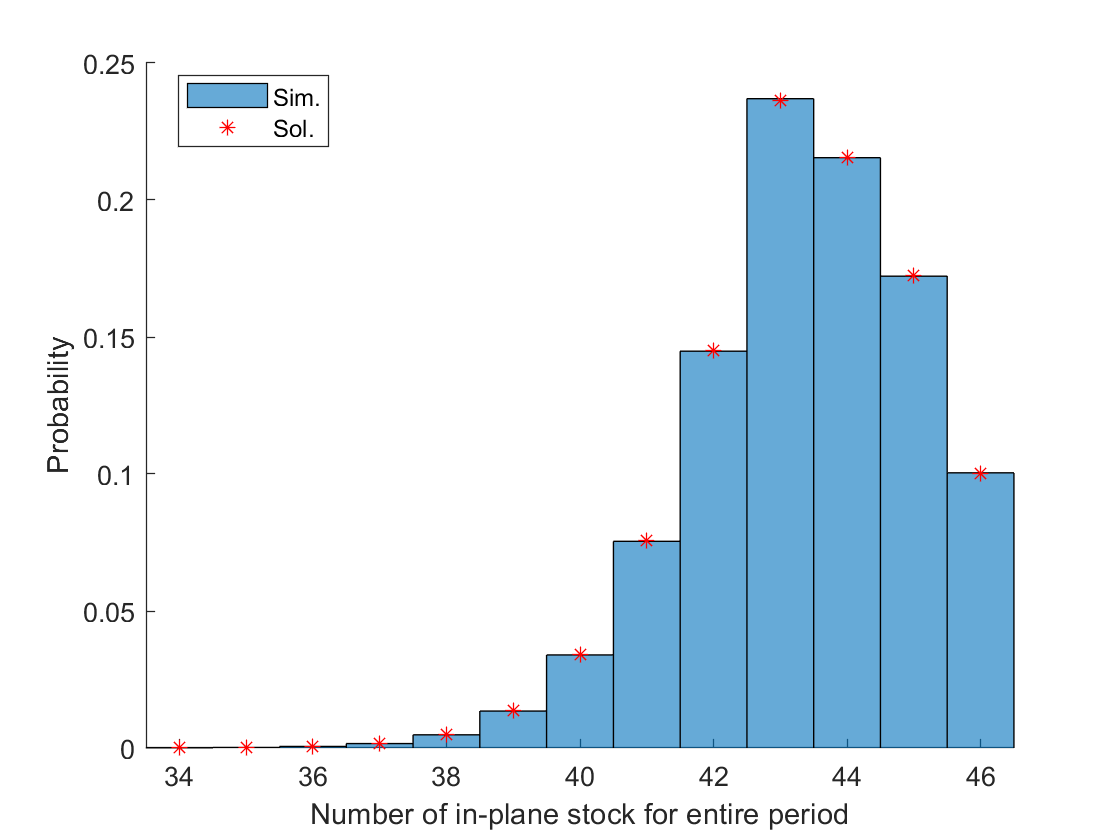}
    \includegraphics[width=.32\textwidth]{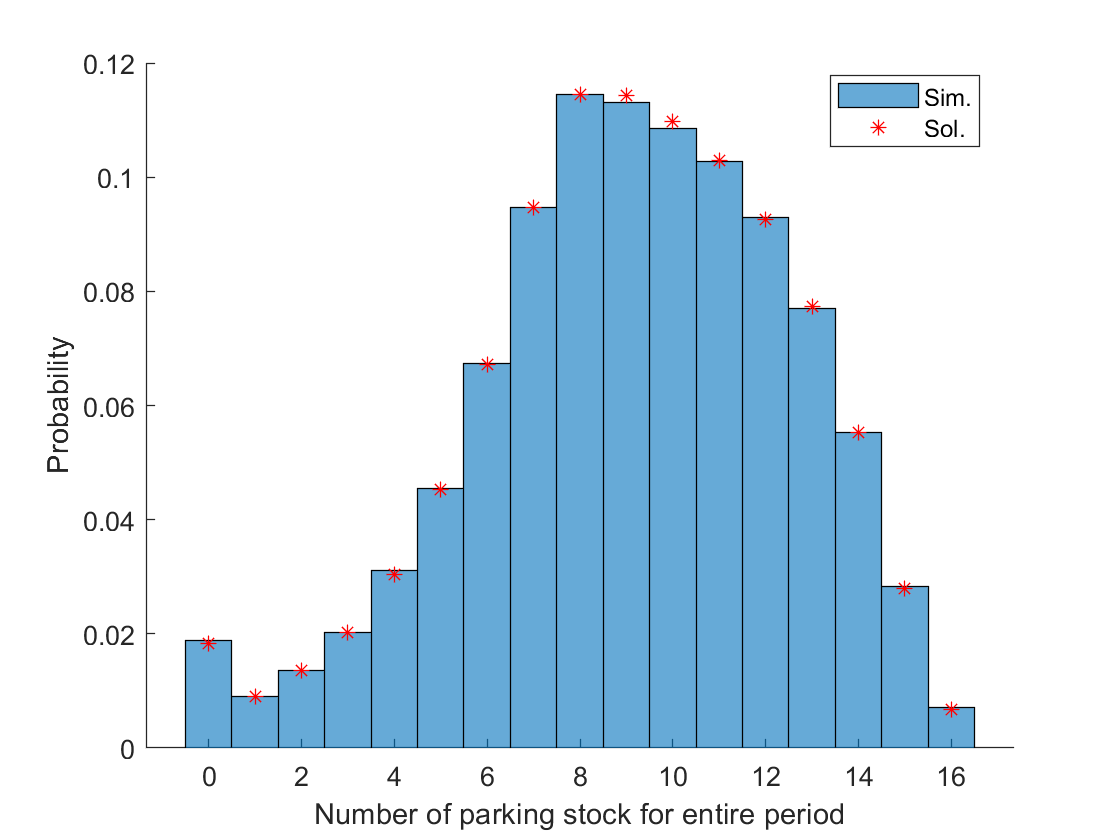}
    \includegraphics[width=.32\textwidth]{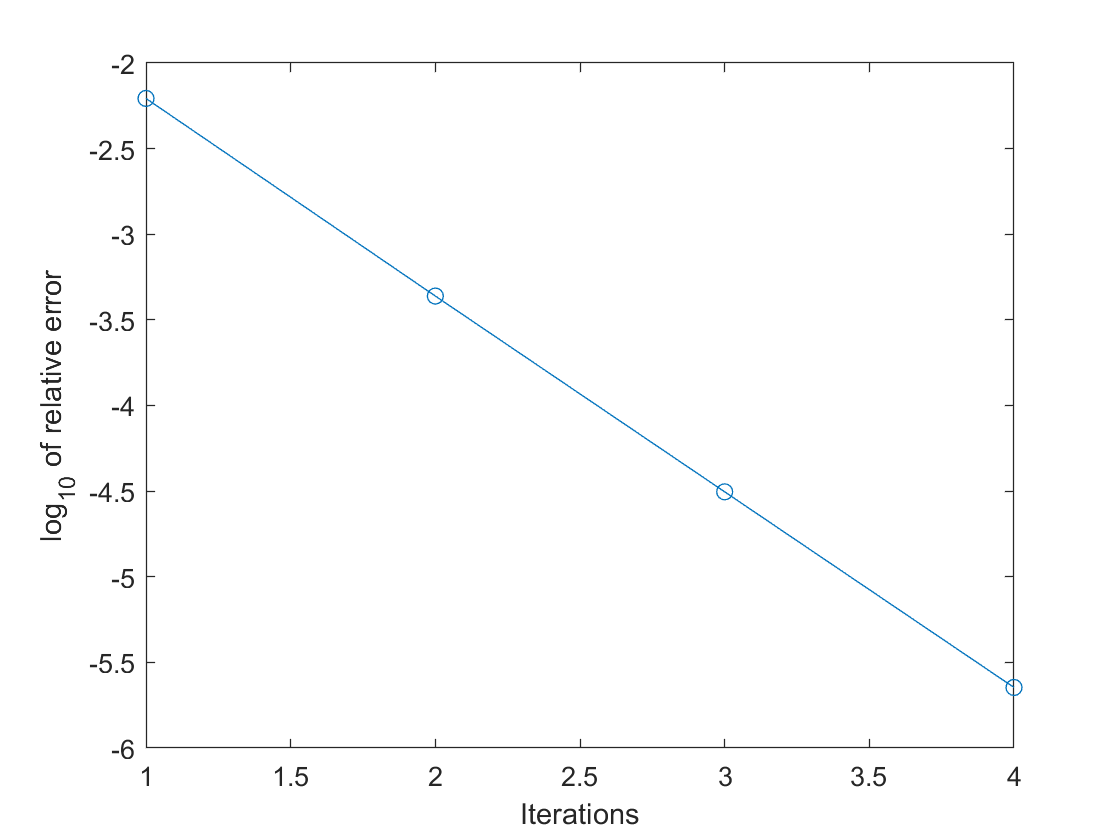}
    \caption{Stationary solution of (a) $\pi^{\text{ir}_i}$ (b) $\pi^{\text{ir}_p}$ and (c) relative error for moderate failure rate}
\end{figure}

\begin{figure}[!h]
    \centering
    \includegraphics[width=.32\textwidth]{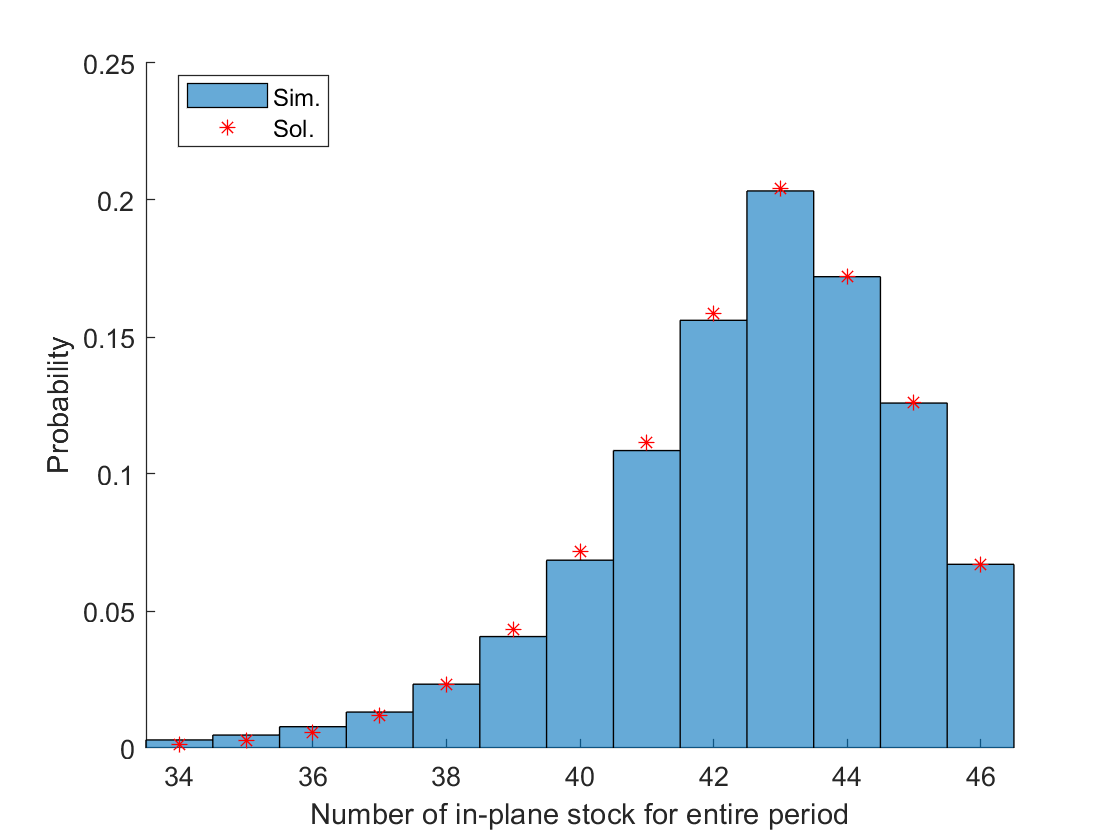}
    \includegraphics[width=.32\textwidth]{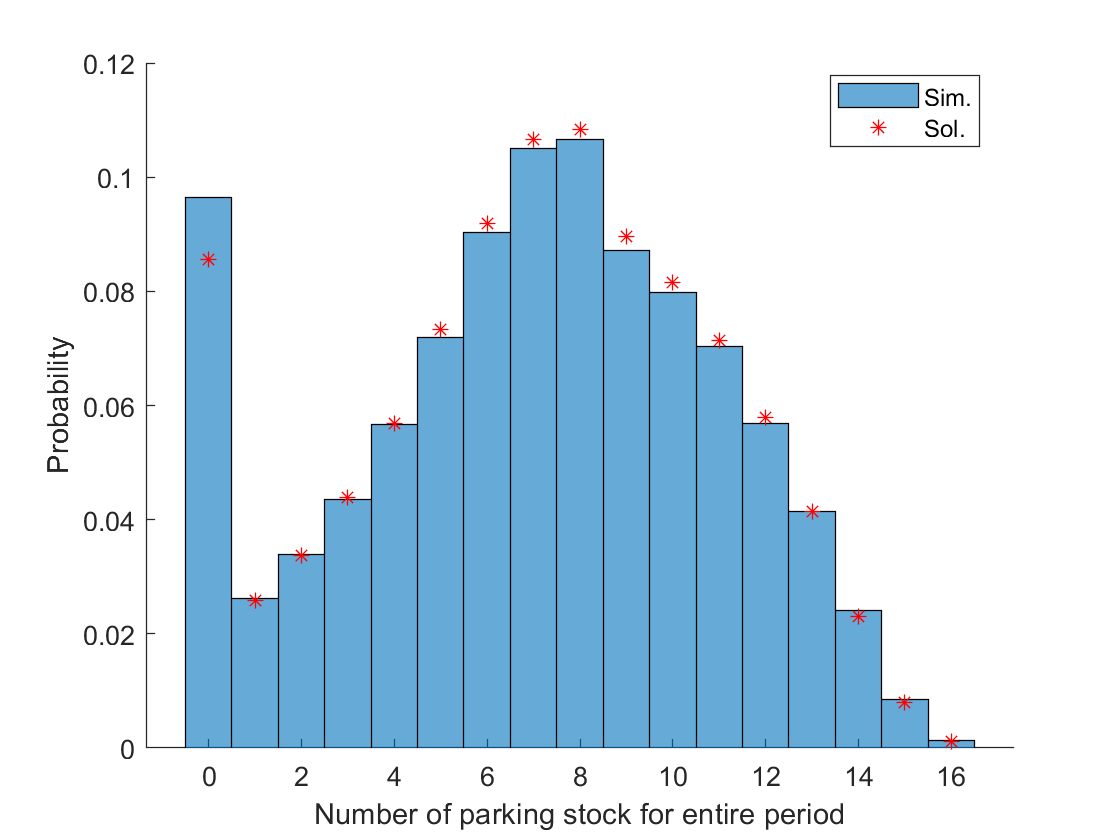}
    \includegraphics[width=.32\textwidth]{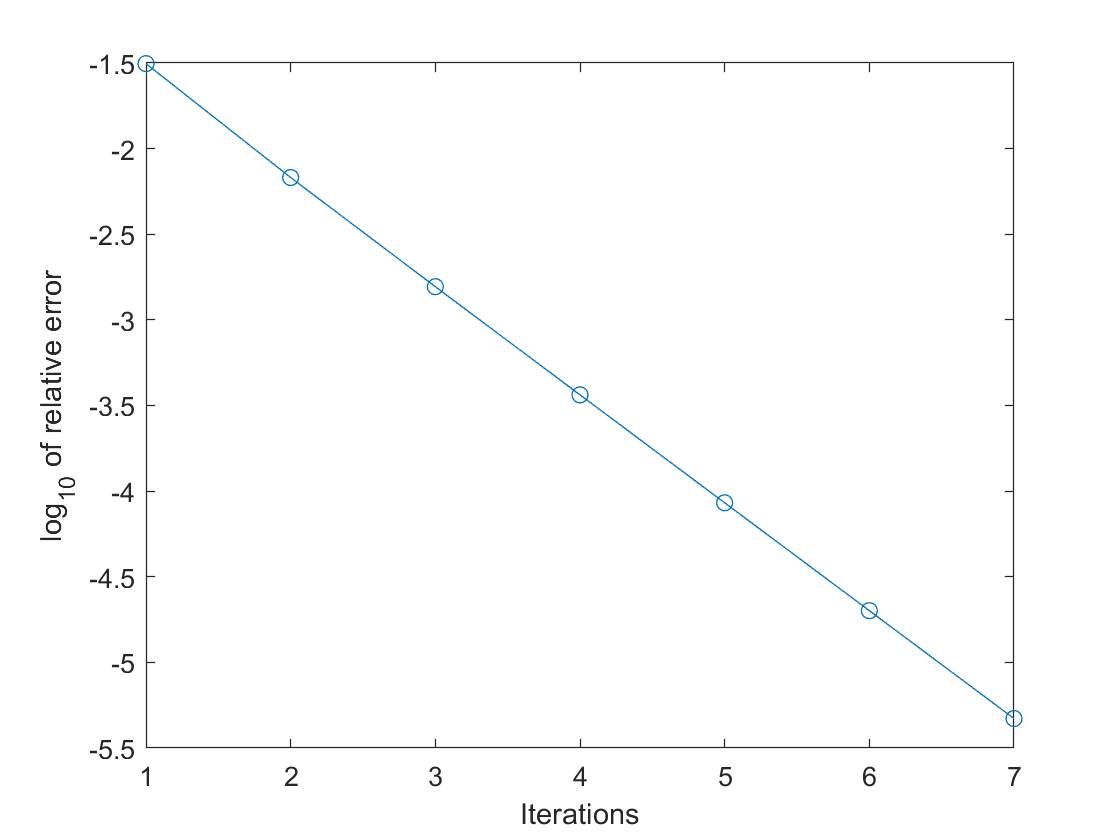}
    \caption{Stationary solution of (a) $\pi^{\text{ir}_i}$ (b) $\pi^{\text{ir}_p}$ and (c) relative error for high failure rate}
\end{figure}
Based on the result figures, the proposed Markov process analysis accurately computes the expected state distribution. However, it loses accuracy as the failure rate increases, or consequently, the parking availability decreases. This is because the assumption of the i.i.d. condition of parking orbits is no longer valid, so the independent analysis loses accuracy. Nevertheless, the proposed method can provide accurate solutions within the practical range of interest. In addition, the number of iterations for convergence is less than 10 for all cases, only requiring a couple of sub-seconds. Therefore, this method can be integrated into the optimization problem without any issues.

\subsection{Simple Optimization of Direct Resupply Strategy}
One representative application of the proposed method is design optimization. We modeled simple costs and constraints for the optimization.
\subsubsection{Cost Modeling}
Let $p_\text{build}$ be the manufacturing price for a single spare satellite. On average, $q$ spares will be sent for every $T_\text{cycle}^\text{dr}$, so the average manufacturing cost for the entire in-plane orbit per unit time is:
\begin{equation}
    c_\text{build} = \frac{1}{T_\text{cycle}^\text{dr}} N_\text{planes} \cdot p_\text{build} \cdot q
\end{equation}
Likewise, every in-plane orbit will have a resupply for every $T_\text{cycle}^\text{dr}$. For launch cost, we assumed that reserving the full payload capacity of an LV will discount the launch cost by $\gamma \%$. Let $p_\text{launch}$ be the launch price for a single satellite and $q_\text{max}$ be the maximum number of spares that can be sent by a single launch. Then, the average launch cost for the entire in-plane orbit per unit time is:
\begin{equation}
    c_\text{launch} = \begin{cases}
        \frac{1}{T_\text{cycle}^\text{dr}} N_\text{planes} \cdot p_\text{launch} \cdot q & \mbox{if } q<q_\text{max} \\
        \frac{1-\gamma}{T_\text{cycle}^\text{dr}} N_\text{planes} \cdot p_\text{launch} \cdot q_\text{max}& \mbox{if } q=q_\text{max} \\
    \end{cases}
\end{equation}
Lastly, having an unnecessarily large number of spare satellites in orbit should be penalized. For this reason, we define the spare holding cost as follows:
\begin{equation}
    c_\text{holding} = N_\text{planes} \cdot \sum_{k=N_\text{sat}+1}^{\bar N_\text{sat}} p_\text{holding} \cdot k \cdot \pi_k^\text{dr}
\end{equation}
where $p_\text{holding}$ is the holding price for a single spare satellite per unit time. Then, the total expected spare strategy cost per unit time becomes:
\begin{equation}
    c_\text{total} = c_\text{build} + c_\text{launch} + c_\text{holding}
\end{equation}

\subsubsection{Constraint Modeling}
To maintain the constellation's designed performance, the average number of operating satellites must be higher than a user-defined threshold. Equivalently, the percentage of time that the number of operating satellites is less than the nominal number of satellites should be smaller than a user-defined threshold $\xi$. This can be modeled as:
\begin{equation}
    P(X < N_\text{sat}) = \sum_{k=0}^{N_\text{sat} - 1} \pi_k^\text{dr} \leq \xi
\end{equation}
In addition, the order size $q$ must satisfy the launch vehicle capacity:
\begin{equation}
    1 \leq q \leq q_\text{max}
\end{equation}
and, though not necessary, we assumed that the reorder level $r$ is larger than or equal to the nominal satellite level.
\begin{equation}
    N_\text{sat} \leq r
\end{equation}

\subsubsection{Optimization Results}
Following parameters are used for optimization of direct resupply policy.
\begin{table}[H]
    \centering
    \begin{tabular}{c|c|c}
    Parameters & Value  & Description \\ \hline
    $p_\text{build}$ & 0.5  & Manufacturing cost per a single spare $[\text{Unit Cost/Sat}]$ \\
    $p_\text{launch}$ & 10 & Launch cost per a single spare $[\text{Unit Cost/Sat}]$ \\
    $p_\text{holding}$ & 0.5 & Launch cost per a single spare per year $[\text{Unit Cost/Sat/Year}]$ \\
    $\xi$ & 0.05 & Performance drop duration requirement $[\%]$ \\
    $\gamma$ & 0.02 & Full launch cost discount factor $[-]$ \\
    $q_\text{max}$ & 6 & The maximum number of spares for a single launch $[\text{\#}]$ \\
    \hline
    \end{tabular}
    \caption{Parameters for the design optimization of the direct resupply method }
\end{table}
\noindent Following graph represents the feasibility and optimality of each design variable $(r,q)$. The color map indicates the total cost of the design set. It turns out that $(r^\ast, q^\ast) = (42,4)$ minimizes the total cost while satisfies the constraints for this scenario.
\begin{figure}[H]
    \centering
    \includegraphics[width=0.38\textwidth]{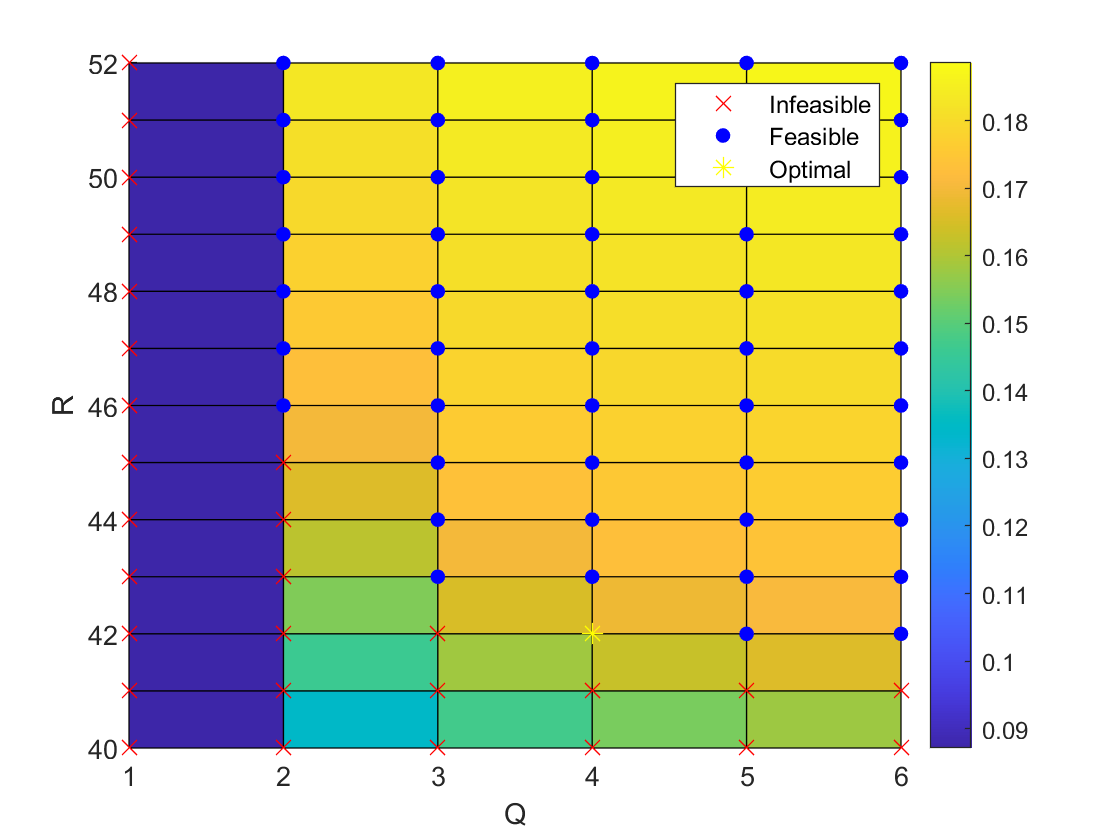}
    \caption{Feasibility and optimality analysis of the direct resupply strategy}
\end{figure}

\section{Conclusion}
In this paper, we have revisited the classical $(r,q)$ inventory control policy and introduced the $(r,q,T)$ inventory control policy. These policies are used to model the direct resupply and indirect resupply strategies. The strategies were mathematically modeled using Markov chains, and their long-term behavior was analyzed using stationary distribution. The proposed method could find an accurate solution very quickly, implying its potential use in design optimization. This paper only covers the simple design optimization of the direct resupply policy.

\section{Acknowledgment}
This research was supported by the Advanced Technology R\&D Center at Mitsubishi Electric Corporation.

\bibliographystyle{AAS_publication}   
\bibliography{references}   

\end{document}